\documentstyle[12pt]{amsart}

\def\opn#1#2{\def#1{\operatorname{#2}}} % to make operators
\opn\chara{char} \opn\length{\ell}
\opn\projdim{proj\,dim} \opn\injdim{inj\,dim} \opn\rank{rank}
\opn\depth{depth} \opn\grade{grade} \opn\height{height}
\opn\embdim{emb\,dim} \opn\codim{codim}

\opn\Tr{Tr} \opn\bigrank{big\,rank}
\opn\superheight{superheight}\opn\lcm{lcm}
\opn\trdeg{tr\,deg}%
\opn\reg{reg} \opn\lreg{lreg}
\opn\div{div} \opn\Div{Div} \opn\cl{cl} \opn\Cl{Cl}
\opn\Spec{Spec} \opn\Supp{Supp} \opn\supp{supp} \opn\Sing{Sing}
\opn\Ass{Ass}
\opn\Ann{Ann} \opn\Rad{Rad} \opn\Soc{Soc}
\opn\Ker{Ker} \opn\Coker{Coker} \opn\Im{Im} \opn\Hom{Hom}
\opn\Tor{Tor} \opn\Ext{Ext} \opn\End{End} \opn\Aut{Aut} \opn\id{id}

\opn\nat{nat}
\opn\pff{pf}%   \pf exists already
\opn\Pf{Pf} \opn\GL{GL} \opn\SL{SL} \opn\mod{mod} \opn\ord{ord}
\opn\aff{aff} \opn\con{conv} \opn\relint{relint} \opn\st{st}
\opn\lk{lk} \opn\cn{cn} \opn\core{core} \opn\vol{vol}
\opn\gr{gr}

\def\pot#1#2{#1[\kern-0.28ex[#2]\kern-0.28ex]}

\opn\dirlim{\underrightarrow{\lim}}
\opn\invlim{\underleftarrow{\lim}}
\def\Implies{\ifmmode\Longrightarrow \else
\unskip${}\Longrightarrow{}$\ignorespaces\fi}
\def\implies{\ifmmode\Rightarrow \else
\unskip${}\Rightarrow{}$\ignorespaces\fi}
\def\iff{\ifmmode\Longleftrightarrow \else
\unskip${}\Longleftrightarrow{}$\ignorespaces\fi}

\let\:=\colon
\newtheorem{Theorem}{Theorem}[section]
\newtheorem{Lemma}[Theorem]{Lemma}
\newtheorem{Corollary}[Theorem]{Corollary}
\newtheorem{Proposition}[Theorem]{Proposition}
\newtheorem{Remark}[Theorem]{Remark}

\newtheorem{Definition}[Theorem]{Definition}

\let\epsilon=\varepsilon
\let\phi=\varphi
\let\kappa=\varkappa
\def\qed{\ifhmode\textqed\fi
\ifmmode\ifinner\quad\qedsymbol\else\dispqed\fi\fi}
\def\textqed{\unskip\nobreak\penalty50
\hskip2em\hbox{}\nobreak\hfil\qedsymbol
\parfillskip=0pt \finalhyphendemerits=0}
\def\dispqed{\rlap{\qquad\qedsymbol}}

\opn\ini{in} \opn\inm{inm} \opn\Sym{Sym}

\begin{document}
\begin{center}
{\large\bf Small slopes of Newton polygon of $L$-function }
\end{center}
\begin{center}
{Fusheng Leng  \ \ \  Banghe Li}
\end{center}
\begin{center}
Academy of Mathematics and Systems Science, Academia Sinica, China
\end{center}

\begin{center}
KLMM
\end{center}

\begin{center}
\bf Abstract
\end{center}

\begin{center}
To understand $L$-function is an important fundamental question in
Number Theory, but there are few specific results on it, especially
the calculation of its Newton polygon. Following Dwork's method it
is hard to calculate an exact example, even on the case of one
variable. There are only three such examples till now, one of which
has some mistakes. In this paper we calculate $L$-functions with
$p$-adic Gauss sums and give a formula in power
series(\textbf{theorem 1.2.}). After that we discuss Newton polygons
$\textrm{NP}(f/\textbf{F}_{p},T)$ of $L$-functions of one variable
polynomials and give a method to calculate its small slopes. We also
obtain the Newton polygon $\textrm{NP}(f/\textbf{F}_{q},T)$ of a
2-variables example with $f=x^{3}+axy+by^{2}$ to illustrate our
method.
\end{center}

\date{}

\section{\bf Introduction}

Let $\textbf{F}_{q}$ be the finite field of $q$ elements with
characteristic $p$ and $\textbf{F}_{q^{k}}$ be the extension of
$\textbf{F}_{q}$ of degree $k$. Let $\zeta_{p}$ be a fixed primitive
$p$-th root of unity in the complex numbers. For any Laurent
polynomial $f(x_{1},\dots,x_{n})\in
\textbf{F}_{q}[x_{1},x_{1}^{-1},\dots,x_{n},x_{n}^{-1}]$, we form
the exponential sum
$$
S_{k}^{*}(f)=\sum_{x_{i}\in\textbf{F}_{q^{k}}^{*}}\zeta_{p}^{\textrm{Tr}_{\textbf{F}_{q^{k}}/\textbf{F}_{p}}(f(x_{1},\cdots,x_{n}))},
\textrm{where} \textbf{
F}_{q^{k}}^{*}=\textbf{F}_{q^{k}}\setminus\{0\}.
$$
The $L$-function is defined by
$$
L^{*}(f,T) =
\textrm{exp}(\sum_{k=1}^{\infty}S_{k}^{*}(f)\frac{T^{k}}{k}) .
$$

To understand the $L$-function is an important fundamental question
in number theory, and since it is very difficult, there are only a
few results on it.

By a theorem of Dwork-Bombieri-Grothendieck,
$$
L^{*}(f,T)
=\frac{\prod_{i=1}^{d_{1}}(1-\alpha_{i}T)}{\prod_{j=1}^{d_{2}}(1-\beta_{j}T)}
$$
is a rational function, where the finitely many numbers $\alpha_{i}$
$(1\leq i\leq d_{1})$ and $\beta_{j}$ $(1\leq j\leq d_{2})$ are
non-zero algebraic integers. Equivalently, for each positive integer
$k$, we have the formula

$$
S_{k}^{*}(f)=\beta_{1}^{k}+\beta_{2}^{k}+ \cdots
+\beta_{d_{2}}^{k}-\alpha_{1}^{k}-\alpha_{2}^{k}-
\cdots-\alpha_{d_{1}}^{k}.
$$

Thus, our fundamental question about the sums $S_{k}^{*}(f)$ is
reduced to understanding the reciprocal zeros $\alpha_{i}$ ($1\leq i
\leq d_{1}$) and the reciprocal poles $\beta_{i}$ ($1\leq j \leq
d_{2}$). When we need to indicate the dependence of the $L$-function
on the ground field $\textbf{F}_{q}$, we will write $L^{*}(f
/\textbf{F}_{q}, T)$.

Without any smoothness condition on $f$, one does not even know
exactly the number $d_{1}$ of zeros and the number $d_{2}$ of poles,
although good upper bounds are available, see [4]. On the other
hand, Deligne's theorem on the Riemann hypothesis [5] gives the
following general information about the nature of the zeros and
poles. For the complex absolute value $\mid \ \mid$, this says

$$
| \alpha_{i}|=q^{\frac{u_{i}}{2}}, | \beta_{j}|=q^{\frac{v_{j}}{2}},
u_{i}\in \textbf{Z}\cap [0, 2n], v_{j}\in \textbf{Z}\cap [0, 2n]
$$

where $\textbf{Z}\cap [0, 2n]$ denotes the set of integers in the
interval $[0, 2n]$. Furthermore, each $\alpha_{i}$(resp. each
$\beta_{j}$) and its Galois conjugates over $\textbf{Q}$ have the
same complex absolute value. For each $l$-adic absolute value $\mid
\ \mid_{l}$ with prime $l\neq p$, the $\alpha_{i}$ and the
$\beta_{j}$ are $l$-adic units:

$$
| \alpha_{i} |_{l}=| \beta_{j} |_{l}=1.
$$

For the remaining prime $p$, it is easy to prove

$$
| \alpha_{i}|_{p}=q^{-r_{i}}, | \beta_{j}|_{p}=q^{-s_{j}}, r_{i}\in
\textbf{Q}\cap [0, n], s_{j}\in \textbf{Q}\cap [0, n].
$$

where we have normalized the $p$-adic absolute value by $| q
|_{p}=q^{-1}$. Deligne's integrality theorem implies the following
improved information:

$$
r_{i}\in \textbf{Q}\cap [0, n], s_{j}\in \textbf{Q}\cap [0, n].
$$

Our fundamental question is then to determine the important
arithmetic invariants $\{u_{i},v_{j},r_{i},s_{j}\}$.

Suppose
$$
f=\sum_{j=1}^{J}a_{j}x^{V_{j}}, a_{j}\neq 0,
$$
where each $V_{j}=(v_{1j}, \cdots, v_{nj})$ is a lattice point in
$\textbf{Z}^{n}$ and the power $x^{V_{j}}$ simply means the product
$x_{1}^{v_{1j}}\cdots x_{n}^{v_{nj}}$. Let $\Delta(f)$ be the convex
closure in $\textbf{R}^{n}$ generated by the origin and the lattice
points $V_{j}$ ($1\leq j\leq J$).

\begin{Definition}
The Laurent polynomial $f$ is called non-degenerate if for each
closed face $\delta$ of $\Delta(f)$ of arbitrary dimension which
does not contain the origin, the $n$ partial derivatives
$$
\{\frac{\partial f^{\delta}}{\partial x_{1}}, \cdots, \frac{\partial
f^{\delta}}{\partial x_{n}}\}
$$
have no common zeros with $x_{1}\cdots x_{n}\neq 0$ over the
algebraic closure of $\textbf{F}_{q}$.
\end{Definition}

If $f$ is non-degenerate, the $L$-function
$L^{*}(f/\textbf{F}_{q},T)^{(-1)^{n-1}}$ is a polynomial of degree
$n!\textbf{V}(f)$ by a theorem of Adolphson-Sperber, where
$\textbf{V}(f)$ denotes the volume of $\Delta(f)$. $^{[2]}$

We are then interested in the Newton polygon of $L$-function. Dwork
gave a method of cohomology theory with $p$-adic on determining
Newton polygon in 1962 and 1964$^{[26][27]}$, and after that
A.Adolphson and S.Sperber developed this method$^{[2]}$. All these
works depend on Dwork's trace formula
$$
S_{k}^{*}(f)=(q^{k}-1)^{n}Tr(\phi^{k})
$$
and the definition of $\phi$ is given by lifting $\textbf{F}_{q}$ to
$\textbf{Q}_{q}$ via a splitting function
$$
\theta(t)=\sum_{m=0}^{\infty}\gamma_{m}t^{m}
$$
where $\phi$ is an endomorphism of some $p$-adic Banach space.

However, there are still few general examples given by the method,
especially when $f$ is a polynomial with one variable. S.Sperber in
1986 gave the Newton polygon of $L^{*}(f,T)$ when $\textrm{deg}
f=3$$^{[16]}$ and fifteen years later S.Hong gave two other examples
with $\textrm{deg} f=4$ and $\textrm{deg} f =6$$^{[9][10]}$, which
still had some mistakes in the last case. $\textrm{deg} f=5$ is more
difficult to calculate than $\textrm{deg} f=6$ in two cases which
need to prove some identical equations by hypergeometric summation
theory, we will show that in our paper.

In 2003, R.Yang calculated the Newton polygon of
$L(f/\textbf{F}_{p},T)$ on a special case when $f=x^{d}+\lambda x$
and $p\equiv -1 \ \textrm{mod}\  d$ for $p$ large enough$^{[22]}$.

In 2004, D.Wan gave a formula of $L$-function with Gauss sum in a
special case when $f$ is diagonal.$^{[28]}$

A Laurent polynomial $f$ is called \textbf{diagonal} if $f$ has
exactly $n$ non-constant terms and $\Delta(f)$ is $n$-dimensional
(necessarily a simplex), then we can write $f(x)$ as
$$
f(x)=\sum_{i=1}^{n}a_{i}x^{V_{i}}, a_{i}\in \textbf{F}_{q}^{*}.
$$

Let $S_{p}$ be the set of $0$ and all rational numbers $a \in (0,1)$
such that $\textrm{ord}_{p}a \geq 0$. For $a, b \in S_{p}$ define
$a+b=c \in S_{p}$ where $c$ is equal to the normal sum $a+b$ $\
\textrm{mod}\  1$. It is not difficult to prove that $(S_{p},+)$
isomorphs to $(\overline{\textbf{F}}_{p},\cdot)$ where
$\overline{\textbf{F}}_{p}$ is the algebraic closure of
$\textbf{F}_{p}$.

Consider the solutions of the following equation
\begin{equation}
(V_{1},\dots,V_{n})(r_{1},\dots,r_{n})^{\textrm{T}}\equiv 0 \
\textrm{mod}\  1, r_{i} \textrm{ rational}, 0\leq r_{i}<1.
\end{equation}

Let $S_{p}(\Delta)$ be the set of solutions $r$ of equation above,
such that $\textrm{ord}_{p}r_{i}\geq 0$ for every $1\leq i\leq n$,
then $S_{p}(\Delta)\subset (S_{p})^{n}$. Let $S_{p}(q,d)$ be the set
of such $r\in S_{p}(\Delta)$ that $(q^{d}-1)r\in \textbf{Z}^{n}$ and
$(q^{d'}-1)r\notin \textbf{Z}^{n}$ for every $1\leq d'< d$. We have
obviously the decomposition
$$
S_{p}(\Delta) = \bigcup_{d\geq 1}S_{p}(q,d).
$$

Let $\chi$ be the Teichm$\ddot{u}$ller character of the
multiplicative group $\textbf{F}_{q}^{*}$. Define Gauss sums over
$\textbf{F}_{q}$ by
$$
G_{k}(q) = -\sum_{a\in
\textbf{F}_{q}^{*}}\chi(a)^{-k}\zeta_{p}^{\textrm{Tr}(a)} (0\leq
k\leq q-2),
$$

D.Wan has proved the following formula when the function $f$ is
diagonal
$$
L^{*}(f/\textbf{F}_{q},T)^{(-1)^{n-1}}=\prod_{d\geq1}\prod_{r\in
S_{p}(q,d)}(1-T^{d}\prod_{i=1}^{n}\chi(a_{i})^{r_{i}(q^{d}-1)}G_{r_{i}(q^{d}-1)}(q^{d}))^{\frac{1}{d}}
.
$$

Note that (1) has finite number solutions since $f$ is diagonal, and
for each of the $d$ points of $r\in S_{p}(q,d)$ the corresponding
factors in this formula are the same. Thus it can be regarded as a
polynomial.

The difficulty on improving this formula to general is, when $f$ is
not diagonal, there will be infinite factors in the formula.

It does not use the "diagonal" condition in the proof of the formula
above and this condition only acts on whether the number of factors
are finite. Note that $S_{p}(q,d)$ is a finite set, does not depend
on whether the Laurent polynomial $f$ is diagonal or not. So we can
get a similar formula in general case, that is
\begin{Theorem}
\label{parameter} Let $f(x)= \sum_{i=1}^{m}a_{i}x^{V_{i}}\in
\textbf{F}_{q}[x_{1},x_{1}^{-1},\dots,x_{n},x_{n}^{-1}]$, where
$a_{i}\neq 0$ for each $i$, and suppose $m>n$, then we have
\begin{equation}
L^{*}(f/\textbf{F}_{q},T)^{(-1)^{n-1}} = \prod_{d\geq 1}\prod_{r\in
S_{p}(q,d)}\prod_{h=0}^{\infty}(1-q^{dh}T^{d}\prod_{i=1}^{m}\chi(a_{i})^{(q^{d}-1)r_{i}}G_{(q^{d}-1)r_{i}}(q^{d}))^{\frac{C_{h+m-n-1}^{m-n-1}}{d}}
\end{equation}
\end{Theorem}

\begin{Remark}
if $r'\equiv q^{s}r \ \textrm{mod}\  1$ for $r,r'\in S_{p}(q,d)$ and
for some integer $s$, the corresponding factors of $r$ and $r'$ in
(2) are the same. Thus, we can remove the power $\frac{1}{d}$ if we
restrict $r$ to run over the $q$-orbits of $S_{p}(q,d)$. Because
$S_{p}(q,d)$ is a finite set, we can easily prove that the right
side of (2) is indeed a power series over $\textbf{Z}_{p}(\pi)$,
where $\textbf{Z}_{p}$ is the ring of $p$-adic integers and $\pi$ be
the unique $(p-1)$-st root of $-p$ satisfying
$\pi\equiv(\zeta_{p}-1) \ \textrm{mod}\  (\zeta_{p}-1)^{2}$.
\end{Remark}

\textbf{Theorem 1.2.} is not the rational function form of
$L$-function, but we will show in this paper that this theorem is
useful on the Newton polygon determination.

Denote the series in \textbf{remark 1.3.} by
$$
L^{*}(f/\textbf{F}_{q},T)^{(-1)^{n-1}}=\sum_{s=0}^{\infty}c_{s}T^{s}.
$$
The main idea in this paper is, to determine the Newton polygon by
calculating $\textrm{ord}_{p}c_{s}$ for every index $s$.

As we have known that
$$
L(f,T)= \textrm{exp}(\sum_{k=1}^{\infty}S_{k}(f)\frac{T^{k}}{k}),
$$
where
$$
S_{k}(f)=\sum_{x_{i}\in\textbf{F}_{q^{k}}}\zeta_{p}^{\textrm{Tr}_{\textbf{F}_{q^{k}}/\textbf{F}_{p}}(f(x_{1},\cdots,x_{n}))},
$$
one will see that $L^{*}(f,T)= (1-T)L(f,T)$ when $f(x)$ is a
polynomial with one variable. Besides, for any $a_{0}\in
\textbf{F}_{q}$, one can easily conclude that $L((f+a_{0}),T)=
L(f,\zeta_{p}^{\textrm{Tr}_{\textbf{F}_{q}/\textbf{F}_{p}}(a_{0})}T)$.
Thus we have $\textrm{NP}((f+a_{0})/\textbf{F}_{q},T)=
\textrm{NP}(f/\textbf{F}_{q},T)$. We can also easily conclude that
such a linear transformation $x+b$ $(b\in \textbf{F}_{q})$ of $x$
does not change the $L$-function, this conclusion will be used to
transform $f$ to reduce our calculation of its Newton polygon.
\begin{center}
\end{center}

In the last of this paper, we will give a new method to calculate
Newton polygons of $L$-functions of one variable polynomials for
first $s$ slopes with $(s-2)(s-1)<2d$ and for every $p$ except some
small values. All of these are based on \textbf{proposition 3.5.}
and \textbf{theorem 6.1.}.
\begin{center}
\end{center}

\smallskip

\section{\bf General theory}
Let $f(x)= \sum_{i=1}^{m}a_{i}x^{V_{i}}\in
\textbf{F}_{q}[x_{1},x_{1}^{-1},\cdots,x_{n},x_{n}^{-1}]$, where
$a_{i}\neq 0$ for each $i$, and suppose $m>n$.

To get the generalization of Wan's formula, we should describe
$S_{p}(q,d)$ first.

Suppose that
\begin{equation}
k_{1}V_{1}+ \cdots + k_{m}V_{m}\equiv 0 \ \textrm{mod}\  (q^{k}-1)
\end{equation}
for a given positive integer $k$ and $0\leq k_{i}\leq q^{k}-2$ for
$i=1,\cdots, m$.

The equation is equivalent to
\begin{equation}
\frac{k_{1}}{q^{k}-1}V_{1}+ \cdots +\frac{k_{m}}{q^{k}-1}V_{m}\equiv
0 \ \textrm{mod}\  1
\end{equation}

Consider the equation
\begin{equation}
r_{1}V_{1}+ \cdots + r_{m}V_{m}\equiv 0 \ \textrm{mod}\  1
\end{equation}
where $r_{i}\in \textbf{Q}$, $0\leq r_{i}<1$ and, if
$r_{i}=\frac{p_{i}}{q_{i}}$ with $(p_{i},q_{i})=1$, then
$(p,q_{i})=1$.

Define $S_{p}(f)$ the solution set of (5). It is clear that if
$(k_{i})$ is a solution of (3), then $(\frac{k_{i}}{q^{k}-1})\in
S_{p}(f)$. Conversely, if $(r_{i})\in S_{p}(f)$, assume
$r_{i}=\frac{p_{i}}{q_{i}}$ with $(p_{i},q_{i})=1$, we will show
that $(r_{i})$ can be written as the form $(\frac{k_{i}}{q^{k}-1})$
for some positive integer $k$.

Since $(p,q_{i})=1$, following the Euler theorem in congruence
theory we have
$$
q^{\lambda_{i}}-1\equiv 0 \ \textrm{mod}\  q_{i}
$$
where $\varphi(q_{i})$ is the Euler function and
$\lambda_{i}|\varphi(q_{i})$ is the smallest positive integer $x$
which satisfies $q^{x}-1\equiv 0 \ \textrm{mod}\  q_{i}$.

Let $\lambda$ be the least common multiple of $\lambda_{1},\cdots,
\lambda_{m}$. Then
$$
q^{\lambda}-1\equiv 0 \ \textrm{mod}\  q_{i}
$$
for all $i$ from 1 to $m$. Thus $r_{i}=\frac{p_{i}}{q_{i}}$ can be
rewritten as the form $\frac{k_{i}}{q^{\lambda}-1}$. That is what we
need.

Let $H_{p}(q,d)$ be the subgroup of $S_{p}(f)$ consisting of all
such
$$
r=(\frac{k_{1}}{q^{d}-1},\cdots, \frac{k_{m}}{q^{d}-1})
$$
with $0\leq k_{i}\leq q^{d}-2$. Then $H_{p}(q,d)\subset H_{p}(q,d')$
if $d| d'$ since $(q^{d}-1)|(q^{d'}-1)$.

Furthermore,
$$
S_{p}(f)=\bigcup_{d\geq 1}H_{p}(q,d).
$$

Define an action
$$
q: r\rightarrow qr=(qr_{1},\cdots, qr_{m}) \ \textrm{mod}\  1
$$
on $S_{p}(f)$. Let $d$ be the number of the elements in the orbit of
$r$ under the action $q$. Then $r\in H_{p}(q,d)$ but $r$ is not in
any $H_{p}(q,d')$ for $d'<d$.

Let
$$
S_{p}(q,d)=H_{p}(q,d)-\bigcup_{d'<d}H_{p}(q,d),
$$
then
$$
H_{p}(q,d)=\bigcup_{d'| d}S_{p}(q,d')
$$
and
$$
S_{p}(f)=\bigcup_{d\geq1}S_{p}(q,d).
$$

It is clear that every subset $S_{p}(q,d)$ is a finite set. In fact
$|S_{p}(q,d)|\leq (q^{d}-1)^{m}$. Furthermore, assume the unique
factorization $d=\prod_{i} p_{i}^{\alpha_{i}}$, following the
principle of inclusion and exclusion we can also prove that
$|S_{p}(q,d)|=(q^{d}-1)-\sum_{i}(q^{\frac{d}{p_{i}}}-1)+\sum_{i\neq
j}(q^{\frac{d}{p_{i}p_{j}}}-1)-\cdots $.

Define Gauss sums over $\textbf{F}_{q}$ by
$$
G_{k}(q) = -\sum_{a\in
\textbf{F}_{q}^{*}}\chi(a)^{-k}\zeta_{p}^{\textrm{Tr}(a)} (0\leq
k\leq q-2).
$$

For each $a\in \textbf{F}_{q}^{*}$, the Gauss sums satisfies the
following interpolation relation
$$
\zeta_{p}^{\textrm{Tr}(a)}=
\sum_{k=0}^{q-2}\frac{G_{k}(q)}{1-q}\chi(a)^{k}.
$$

To get the generalization of Wan's formula, we also need a formula
on Gauss sums. That is

\begin{Theorem}
(Hasse-Davenport) For every positive integer $k$,
$$
G_{r(q^{dk}-1)}(q^{dk})= G_{r(q^{d}-1)}(q^{d})^{k}.
$$
\end{Theorem}

Then we can calculate that
\begin{eqnarray*}
S_{1}^{*}(f)&=& \sum_{x_{j}\in
\textbf{F}_{q}^{*}}\zeta_{p}^{\textrm{Tr}(f(x))}\\
&=& \sum_{x_{j}\in
\textbf{F}_{q}^{*}}\prod_{i=1}^{m}\zeta_{p}^{\textrm{Tr}(a_{i}x^{V_{i}})}\\
&=& \sum_{x_{j}\in
\textbf{F}_{q}^{*}}\prod_{i=1}^{m}\sum_{k_{i}=0}^{q-2}
\frac{G_{k_{i}}(q)}{1-q}\chi(a_{i})^{k_{i}}\chi(x^{V_{i}})^{k_{i}}\\
&=&\sum_{k_{1}=0}^{q-2}\cdots\sum_{k_{m}=0}^{q-2}(\prod_{i=1}^{m}\frac{G_{k_{i}}(q)}{1-q}\chi(a_{i})^{k_{i}})\sum_{x_{j}\in
\textbf{F}_{q}^{*}}\chi(x^{k_{1}V_{1}+\cdots+k_{m}V_{m}}).
\end{eqnarray*}

Note that $\sum_{x_{j}\in
\textbf{F}_{q}^{*}}\chi(x^{k_{1}V_{1}+\cdots+k_{m}V_{m}})= 0$ unless
$k_{1}V_{1}+\cdots+k_{m}V_{m}\equiv 0 \ \textrm{mod}\  q-1$, if this
condition holds, the value of $\sum_{x_{j}\in
\textbf{F}_{q}^{*}}\chi(x^{k_{1}V_{1}+\cdots+k_{m}V_{m}})$ will be
$(q-1)^{n}$. Thus,
$$
S_{1}^{*}(f)=
(-1)^{n}(1-q)^{n-m}\sum_{k_{1}V_{1}+\cdots+k_{m}V_{m}\equiv 0 \
\textrm{mod}\  q-1}\prod_{i=1}^{m}\chi(a_{i})^{k_{i}}G_{k_{i}}(q).
$$
Replacing $q$ by $q^{k}$, one gets a formula for the exponential sum
$S_{k}^{*}(f)$ over the $k$-th extension field $\textbf{F}_{q^{k}}$:

$$
S_{k}^{*}(f)= (-1)^{n}(1-q^{k})^{n-m}\sum_{r\in
H_{p}(q,k)}\prod_{i=1}^{m}\chi(a_{i})^{r_{i}(q^{k}-1)}G_{r_{i}(q^{k}-1)}(q^{k})
$$
\begin{equation}
= \sum_{k'| k}\sum_{r\in
S_{p}(q,k')}(-1)^{n}(1-q^{k})^{n-m}\prod_{i=1}^{m}\chi(a_{i})^{r_{i}(q^{k}-1)}G_{r_{i}(q^{k}-1)}(q^{k}).
\end{equation}
\begin{center}
\end{center}

Since
$$
L^{*}(f/\textbf{F}_{q},T) =
\textrm{exp}(\sum_{k=1}^{\infty}S_{k}^{*}(f)\frac{T^{k}}{k}),
$$
by (6) the equation above is equal to
$$
\prod_{d\geq 1}\prod_{r\in
S_{p}(q,d)}\textrm{exp}(\sum_{k=1}^{\infty}\frac{T^{kd}}{kd}(-1)^{n}(1-q^{kd})^{n-m}
\prod_{i=1}^{m}\chi(a_{i})^{r_{i}(q^{kd}-1)}G_{r_{i}(q^{kd}-1)}(q^{kd})).
$$
Following the Hasse-Davenport relation we rewrite it by
\begin{eqnarray*}
L^{*}(f/\textbf{F}_{q},T)&=& \prod_{d\geq 1}\prod_{r\in
S_{p}(q,d)}\textrm{exp}(\sum_{k=1}^{\infty}\frac{T^{kd}}{kd}(-1)^{n}(1-q^{kd})^{n-m}
\prod_{i=1}^{m}\chi(a_{i})^{kr_{i}(q^{d}-1)}G_{r_{i}(q^{d}-1)}(q^{d})^{k})\\
&=& \prod_{d\geq 1}\prod_{r\in
S_{p}(q,d)}\textrm{exp}(\sum_{k=1}^{\infty}\frac{T^{kd}}{kd}(-1)^{n}\sum_{h=0}^{\infty}
C_{h+m-n-1}^{m-n-1}q^{kdh}\prod_{i=1}^{m}\chi(a_{i})^{kr_{i}(q^{d}-1)}G_{r_{i}(q^{d}-1)}(q^{d})^{k})\\
&=& \prod_{d\geq 1}\prod_{r\in
S_{p}(q,d)}\prod_{h=0}^{\infty}\textrm{exp}(\sum_{k=1}^{\infty}\frac{T^{kd}}
{kd}(-1)^{n}C_{h+m-n-1}^{m-n-1}q^{kdh}\prod_{i=1}^{m}\chi(a_{i})^{kr_{i}(q^{d}-1)}G_{r_{i}(q^{d}-1)}(q^{d})^{k})\\
&=& (\prod_{d\geq 1}\prod_{r\in
S_{p}(q,d)}\prod_{h=0}^{\infty}(1-q^{dh}T^{d}\prod_{i=1}^{m}\chi(a_{i})^{r_{i}(q^{d}-1)}G_{r_{i}(q^{d}-1)}(q^{d}))^{\frac{C_{h+m-n-1}^{m-n-1}}{d}})^{(-1)^{n-1}}.
\end{eqnarray*}

This is the proof of \textbf{Theorem 1.2.}.
\begin{center}
\end{center}

We denote the $h=0$-part in (2) as
$$
L_{0}^{*}(f/\textbf{F}_{q},T)=\prod_{d\geq 1}\prod_{r\in
S_{p}(q,d)}(1-T^{d}\prod_{i=1}^{m}\chi(a_{i})^{(q^{d}-1)r_{i}}G_{(q^{d}-1)r_{i}}(q^{d}))^{\frac{1}{d}}.
$$

Recall that $ S_{1}^{*}(f)=
(-1)^{n}(1-q)^{n-m}\sum_{k_{1}V_{1}+\cdots+k_{m}V_{m}\equiv 0 \
\textrm{mod}\  q-1}\prod_{i=1}^{m}\chi(a_{i})^{k_{i}}G_{k_{i}}(q)$
the $h=0$-part in (2) is indeed the
$\textrm{exp}(\sum_{k=1}^{\infty}S_{k}^{*}(f)\frac{T^{k}}{k})$
replacing $S_{k}^{*}(f)$ by $(1-q^{k})^{m-n}S_{k}^{*}(f)$.

\begin{Theorem}
Suppose $f$ is non-degenerate, then the Newton polygon of
$L^{*}(f/\textbf{F}_{q},T)$ is the same as
$L_{0}^{*}(f/\textbf{F}_{q},T)$ up to the slopes smaller than 1.
Especially, if $f(x)= \sum_{i=1}^{m}a_{i}x^{d_{i}}$, where
$0<d_{1}<\cdots<d_{m}= d$, $a_{i}\in\textbf{F}_{q}^{*}$ and $d\neq
0\ \textrm{mod}\  p$, then $f$ is non-degenerate and the Newton
polygon of $L^{*}(f/\textbf{F}_{q},T)$ is the same as
$L_{0}^{*}(f/\textbf{F}_{q},T)$ up to $d$.
\end{Theorem}
\begin{pf}
Since

$$
L_{0}^{*}(f/\textbf{F}_{q},T)=\textrm{exp}(\sum_{k=1}^{\infty}(1-q^{k})^{m-n}S_{k}^{*}(f)\frac{T^{k}}{k}),
$$
we have

\begin{eqnarray*}
L_{0}^{*}(f/\textbf{F}_{q},T)
&=&\textrm{exp}(\sum_{k=1}^{\infty}\sum_{i=0}^{m-n}C_{m-n}^{i}(-1)^{i}q^{ki}\frac{S_{k}^{*}(f)T^{k}}
{k})\\
&=&\prod_{i=0}^{m-n}(\textrm{exp}(\sum_{k=1}^{\infty}\frac{S_{k}^{*}(f)(Tq^{i})^{k}}
{k}))^{C_{m-n}^{i}(-1)^{i}}\\
&=&\frac{L^{*}(f/\textbf{F}_{q},T)L^{*}(f/\textbf{F}_{q},q^{2}T)^{C_{m-n}^{2}}\cdots}{L^{*}(f/\textbf{F}_{q},qT)^{m-n}L^{*}(f/\textbf{F}_{q},q^{3}T)^{C_{m-n}^{3}}\cdots}.
\end{eqnarray*}

Recall that

$$
L^{*}(f,T)
=\frac{\prod_{i=1}^{d_{1}}(1-\alpha_{i}T)}{\prod_{j=1}^{d_{2}}(1-\beta_{j}T)}
$$

where

$$
\textrm{ord}_{q}\alpha_{i}=r_{i}, \textrm{ord}_{q}\beta_{j}=s_{j},
r_{i}\in \textbf{Q}\cap [0, n], s_{j}\in \textbf{Q}\cap [0, n].
$$

with normalized the $p$-adic order by $\textrm{ord}_{q}q=1$, then
the reciprocal zeros $q^{k}\alpha_{i}$ ($1\leq i \leq d_{1}$) and
the reciprocal poles $q^{k}\beta_{i}$ ($1\leq j \leq d_{2}$) of
$L^{*}(f/\textbf{F}_{q},q^{k}T)$ have the $p$-adic orders

$$
\textrm{ord}_{q}(q^{k}\alpha_{i})=r_{i}+k,
\textrm{ord}_{q}(q^{k}\beta_{j})=s_{j}+k, r_{i}\in \textbf{Q}\cap
[0, n], s_{j}\in \textbf{Q}\cap [0, n].
$$

Thus, the only reciprocal zeros or reciprocal poles in
$L_{0}^{*}(f/\textbf{F}_{q},T)$ for which the $p$-adic orders
smaller than 1 must appear in $L^{*}(f/\textbf{F}_{q},T)$ of the
right side of the equation above.

Then we obtain the theorem.
\end{pf}

Following \textbf{theorem 2.2.}, to calculate the Newton polygon of
$L^{*}(f/\textbf{F}_{q},T)$ where $f(x)=
\sum_{i=1}^{m}a_{i}x^{d_{i}}$ we should only determine the
$\textrm{ord}_{q}$-value of $L_{0}^{*}(f/\textbf{F}_{q},T)$'s first
$d$-terms.

\section{\bf Presentation of $L_{0}^{*}(f/\textbf{F}_{q},T)$'s term}

To express clearly, we first give some symbols that will be used
below.

\begin{Definition}
For an arbitrary real number $x$, define $\{x\}$ satisfying that
$0\leq \{x\}<1$ and $\{x\}\equiv x \ \textrm{mod}\  1$.
\end{Definition}

\begin{Definition}
Any nonnegative integer $k$ can be written as the form
$$
k=k_{0}+pk_{1}+ \cdots +p^{l-1}k_{l-1}
$$
uniquely, where $k_{t}$ is an integer and $0\leq k_{t}\leq p-1$ for
each $0\leq t<l$ and $k_{l-1}\neq 0$. Define a function $\sigma$ on
the nonnegative integer set to itself such that
$$
\sigma(k)=\sum_{t=0}^{l-1}k_{t}.
$$
\end{Definition}
\begin{center}
\end{center}

Following these symbols we can express the Gross-Koblitz formula as
below:

\begin{Theorem}
(Gross-Koblitz) Suppose $p\geq 2$ prime and $q=p^{a}$. Let $\pi$ be
the unique $(p-1)$-st root of $-p$ satisfying
$$
\pi\equiv(\zeta_{p}-1) \ \textrm{mod}\  (\zeta_{p}-1)^{2}.
$$
Then
$$
G_{k}(q)=\pi^{\sigma(k)}\Pi_{j=0}^{a-1}\Gamma_{p}(\{\frac{p^{j}k}{q-1}\}).
$$
\end{Theorem}
\begin{center}
\end{center}

Let $f(x)= \sum_{i=1}^{m}a_{i}x^{d_{i}}$, where
$0<d_{1}<\cdots<d_{m}= d$, $\textrm{gcd}(p, d)=1$,
$a_{i}\in\textbf{F}_{q}^{*}$. Thus $f$ is non-degenerate and
$L^{*}(f/\textbf{F}_{q},T)$ is a polynomial of degree $d$.

Recall the definition of $L_{0}^{*}(f/\textbf{F}_{p},T)$ we denote
$L_{0}^{*}(f/\textbf{F}_{p},T)=\sum_{s=0}^{\infty}c_{s}T^{s}$. For a
given integer $s$, the $c_{s}$ can be expressed as the sum of all
these terms:
\begin{equation}
\prod_{\sum s_{j}=s,r_{j}\in
S_{p}(p,s_{j})}(-\prod_{i=1}^{m}\chi(a_{i})^{(p^{s_{j}}-1)r_{ji}}G_{(p^{s_{j}}-1)r_{ji}}(p^{s_{j}})).
\end{equation}

Following Gross-Koblitz formula (7) can be written as
$$
\prod_{\sum s_{j}=s,r_{j}\in
S_{p}(p,s_{j})}(-\prod_{i=1}^{m}\chi(a_{i})^{(p^{s_{j}}-1)r_{ji}}G_{(p^{s_{j}}-1)r_{ji}}(p^{s_{j}}))
$$
\begin{equation}
=\prod_{i=1}^{m}(\chi(a_{i})^{\sum_{j}\sigma((p^{s_{j}}-1)r_{ji})}\pi^{\sum_{j}\sigma((p^{s_{j}}-1)r_{ji})})
\cdot\prod_{j}(-\prod_{i=1}^{m}\prod_{t=0}^{s_{j}-1}\Gamma_{p}(\{p^{t}r_{ji}\}))
.
\end{equation}

There are infinite number of such addends, but only finite number of
them have the $\textrm{ord}_{p}$-value smaller than a given number.

Consider (8), let
$$
k_{ji}=(p^{s_{j}}-1)r_{ji}=\sum_{t=0}^{s_{j}-1}k_{ji}[t]p^{t}
$$
and
$$
\sum_{i=1}^{m}d_{i}k_{ji}[t]=u_{j}[t]p-v_{j}[t]
$$
where $0\leq k_{ji}[t]\leq p-1$ and $0\leq v_{j}[t]\leq p-1$ and
$0\leq u_{j}[t]\leq \sum_{i=1}^{m}d_{i}$.

\begin{Proposition}
Suppose $p\geq \sum_{i=1}^{m}d_{i}$. For any index $j$,
\begin{equation}
u_{j}[t-1]=v_{j}[t]
\end{equation}
for every $0\leq t\leq s_{j}-1$.

Conversely, let
$$
k_{ji}=\sum_{t=0}^{s_{j}-1}k_{ji}[t]p^{t}
$$
and
$$
\sum_{i=1}^{m}d_{i}k_{ji}[t]=u_{j}[t]p-v_{j}[t]
$$
for $0\leq k_{ji}[t]\leq p-1$ and $0\leq v_{j}[t]\leq p-1$ and
$0\leq u_{j}[t]\leq \sum_{i=1}^{m}d_{i}$. if the condition
$$
u_{j}[t-1]=v_{j}[t]
$$
for every $0\leq t\leq s_{j}-1$ is achieved, then
$$
r_{j}=(\frac{k_{j1}}{p^{s_{j}}-1}, \cdots,
\frac{k_{jm}}{p^{s_{j}}-1})\in S_{p}(p,s_{j}).
$$
\end{Proposition}

\begin{pf}
Since $r_{j}\in S_{p}(p,s_{j})$ we have
$$
\sum_{i=1}^{m}d_{i}k_{ji}\equiv 0 \ \textrm{mod}\  p^{s_{j}}-1.
$$
On the other hand,
\begin{eqnarray*}
\sum_{i=1}^{m}d_{i}k_{ji}
&=&\sum_{i=1}^{m}\sum_{t=0}^{s_{j}-1}d_{i}k_{ji}[t]p^{t}\\
&=&\sum_{t=0}^{s_{j}-1}\sum_{i=1}^{m}d_{i}k_{ji}[t]p^{t}\\
&=&\sum_{t=0}^{s_{j}-1}(u_{j}[t]p-v_{j}[t])p^{t}\\
&=&u_{j}[s_{j}-1](p^{s_{j}}-1)+\sum_{t=0}^{s_{j}-1}(u_{j}[t-1]-v_{j}[t])p^{t},
\end{eqnarray*}
where $u_{j}[-1]=u_{j}[s_{j}-1]$.

Then
$$
\sum_{i=1}^{m}d_{i}k_{ji}\equiv 0 \ \textrm{mod}\  p^{s_{j}}-1
$$
is equivalent to
\begin{equation}
\sum_{t=0}^{s_{j}-1}(u_{j}[t-1]-v_{j}[t])p^{t}\equiv 0 \
\textrm{mod}\ p^{s_{j}}-1.
\end{equation}

Recall that
$$
\sum_{i=1}^{m}d_{i}k_{ji}[t]=u_{j}[t]p-v_{j}[t],
$$
by $p\geq \sum_{i=1}^{m}d_{i}$, $0\leq k_{ji}[t]\leq p-1$ and $0\leq
v_{j}[t]\leq p-1$ we have
\begin{eqnarray*}
u_{j}[t]
&=&\frac{1}{p}(\sum_{i=1}^{m}d_{i}k_{ji}[t]+v_{j}[t])\\
&\leq&\frac{1}{p}(\sum_{i=1}^{m}d_{i}k_{ji}[t]+(p-1))\\
&\leq&\frac{1}{p}(1+\sum_{i=1}^{m}d_{i})(p-1)\\
&\leq&\frac{1}{p}(1+p)(p-1)\\
&<&p
\end{eqnarray*}
i.e.
$$
u_{j}[t]\leq p-1.
$$

Then we have
\begin{eqnarray*}
\sum_{t=0}^{s_{j}-1}(u_{j}[t-1]-v_{j}[t])p^{t}
&\leq&\sum_{t=0}^{s_{j}-1}u_{j}[t-1]p^{t}\\
&\leq&\sum_{t=0}^{s_{j}-1}(p-1)p^{t}\\
&=&p^{s_{j}}-1
\end{eqnarray*}

If the equality is achieved, it obtains that
\begin{equation}
v_{j}[t]=0, u_{j}[t]=p-1
\end{equation}
for every $0< t\leq s_{j}-1$. Furthermore, by
$\sum_{i=1}^{m}d_{i}k_{ji}[t]=u_{j}[t]p-v_{j}[t]$ and $0\leq
k_{ji}[t]\leq p-1$, (11) means that
$$
p=\sum_{i=1}^{m}d_{i}
$$
and
$$
k_{ji}[t]=p-1
$$
for every $0< t\leq s_{j}-1$ and $1\leq i\leq m$, but it is
contradictory to $k_{ji}<p^{s_{j}}-1$.

Thus
\begin{equation}
\sum_{t=0}^{s_{j}-1}(u_{j}[t-1]-v_{j}[t])p^{t}<p^{s_{j}}-1.
\end{equation}

Besides,
\begin{eqnarray*}
\sum_{t=0}^{s_{j}-1}(u_{j}[t-1]-v_{j}[t])p^{t}
&\geq&-\sum_{t=0}^{s_{j}-1}v_{j}[t]p^{t}\\
&\geq&-\sum_{t=0}^{s_{j}-1}(p-1)p^{t}\\
&=&-(p^{s_{j}}-1)
\end{eqnarray*}
The equality can not be achieved since
$u_{j}[t]p-v_{j}[t]=\sum_{i=1}^{m}d_{i}k_{ji}[t]\geq 0$ for every
$0< t\leq s_{j}-1$, i.e.
\begin{equation}
-(p^{s_{j}}-1)<\sum_{t=0}^{s_{j}-1}(u_{j}[t-1]-v_{j}[t])p^{t}.
\end{equation}

Following (10), (12) and (13) we obtain that
$$
\sum_{t=0}^{s_{j}-1}(u_{j}[t-1]-v_{j}[t])p^{t}=0,
$$
i.e.
$$
u_{j}[t-1]=v_{j}[t]
$$
for every $0\leq t\leq s_{j}-1$.

To prove the remainder of this Proposition, we should only show that
$$
\sum_{i=1}^{m}d_{i}k_{ji}\equiv 0 \ \textrm{mod}\  p^{s_{j}}-1
$$
when the condition is achieved. This is obvious.
\end{pf}
\begin{center}
\end{center}

Following the notations above we construct a table with many blocks
as below:
\begin{center}
\end{center}
$\left(
\begin{array}{cccc}
k_{11}[0] & k_{11}[1] & \cdots & k_{11}[s_{1}-1] \\
k_{12}[0] & k_{12}[1] & \cdots & k_{12}[s_{1}-1] \\
\cdots & \cdots & \cdots & \cdots \\
k_{1m}[0] & k_{1m}[1] & \cdots & k_{1m}[s_{1}-1] \\
\end{array}
\right)$ $\cdots $ $\left(
\begin{array}{cccc}
k_{j1}[0] & k_{j1}[1] & \cdots & k_{j1}[s_{j}-1] \\
k_{j2}[0] & k_{j2}[1] & \cdots & k_{j2}[s_{j}-1] \\
\cdots & \cdots & \cdots & \cdots \\
k_{jm}[0] & k_{jm}[1] & \cdots & k_{jm}[s_{j}-1] \\
\end{array}
\right)$ $\cdots $
\begin{center}
\end{center}
where $\sum_{j}s_{j}=s$.
\begin{center}
\end{center}

A block
\begin{center}
\end{center}
$\left(
\begin{array}{cccc}
k_{j1}[0] & k_{j1}[1] & \cdots & k_{j1}[s_{j}-1] \\
k_{j2}[0] & k_{j2}[1] & \cdots & k_{j2}[s_{j}-1] \\
\cdots & \cdots & \cdots & \cdots \\
k_{jm}[0] & k_{jm}[1] & \cdots & k_{jm}[s_{j}-1] \\
\end{array}
\right)$
\begin{center}
\end{center}
in the table is corresponding to such a factor
$$
-\prod_{i=1}^{m}\chi(a_{i})^{(p^{s_{j}}-1)r_{ji}}G_{(p^{s_{j}}-1)r_{ji}}(p^{s_{j}})
$$
in (7). This table is determined by the $c_{s}$' term which we have
chosen.

Conversely, if such a table is given, satisfying the relation in
\textbf{Proposition 3.4.}, each block is corresponding to a
$r_{j}\in S_{p}(p,s_{j})$, where $s_{j}$ is the number of the
block's columns, and each of such $r_{j}$ is in the different
orbits, then the table determines a term of $c_{s}$.

Following (8) a term of $c_{s}$ corresponding to such a table above
has the $\textrm{ord}_{p}$-value equal to
$$
\textrm{ord}_{p}\pi\cdot (\sum_{j}\sigma((p^{s_{j}}-1)r_{ji})),
$$
which is $\frac{1}{p-1}$ products the sum of all $k_{ji}[t]$ that
appear in the table, i.e.
$$
\frac{1}{p-1}\sum_{j}\sum_{i=1}^{m}\sum_{t=0}^{s_{j}-1}k_{ji}[t].
$$

Let $T(c_{s})$ be the set consisting of all the terms of $c_{s}$,
and $I(c_{s})=\{\textrm{ord}_{p}c\mid c\in T(c_{s})\}$. Rewrite
$I(c_{s})$ to $\{\alpha_{1},\cdots, \alpha_{k},\cdots\}$, where
$\alpha_{k}<\alpha_{k'}$ when $k<k'$, and denote $T_{k}=\{c\in
T(c_{s})\mid\textrm{ord}_{p}c=\alpha_{k}\}$. To calculate
$\textrm{ord}_{p}c_{s}$, we should calculate
$$
\textrm{ord}_{p}\sum_{c\in T_{1}}c.
$$
If it is equal to $\alpha_{1}$, then $\textrm{ord}_{p}c_{s}$ is also
equal to $\alpha_{1}$, otherwise we will show in the next two
sections that
$$
\textrm{ord}_{p}\sum_{c\in T_{1}}c\geq \alpha_{1}+1,
$$
and we should calculate
$$
\textrm{ord}_{p}\sum_{c\in T_{2}}c,
$$
and so on.

We will point out that the sum of all terms as form (8) whose
$\textrm{ord}_{p}$-value smaller than a certain number has the
$\textrm{ord}_{p}$-value greater than 1, that means we can submit
those terms when calculating $\textrm{ord}_{p}c_{s}$ and begin from
some greater value of $\alpha_{k}$. Indeed, we can get

\begin{Proposition}
Suppose $p\geq \sum_{i=1}^{m}d_{i}$. For $s>0$, we have
$$
\textrm{ord}_{p}\sum c \geq 1+\frac{s-1}{d}
$$
for all such $c\in T(c_{s})$ that, there are two same $u$-values
$u_{j}[t]=u_{j'}[t']$ in the table corresponding to $c$.

Furthermore, if  $(s-2)(s-1)<2d$, to calculate
$\textrm{ord}_{p}c_{s}$ we should only begin from
$\textrm{ord}_{p}\sum c$ for all such $c\in T(c_{s})$ that
corresponds to $s$ different $u$-values as $0, 1, \cdots, s-1$
respectively.
\end{Proposition}

To prove this proposition, we will introduce some definitions and
conclusions in the next two sections.
\begin{center}
\end{center}

Denote $w_{1},\dots , w_{s} $ the $s$ column vectors of such
$(k_{ji}[t])_{i=1,\dots,m}$. It is only a symbol in general that
have no $k_{ji}[t]$'s values yet. When the values of $k_{ji}[t]$ are
fixed, we denote $s$ vectors $\overline{w_{1}},\dots ,
\overline{w_{s}} $ as the values of $w_{1},\dots , w_{s} $
respectively. Note that every $w_{j}$ is different since it is just
a symbol, and some of $\overline{w_{j}}$ may have the same value of
$u$ or $v$ and, even they are all the same as vectors.
\begin{center}
\end{center}

\section{\bf $f$-simple permutation on the symmetric group}

In this section, $f$ is a given map defined on the set $\{1, 2, 3,
\dots, n\}$ and maps to an arbitrary set, and $S_{n}$ is the
symmetric group of $\{1, 2, 3, \dots, n\}$. Note that a permutation
can be written as the product of several separated \textbf{cycles}
uniquely.

It is easily to see that all $\sigma \in S_{n}$ satisfying
$$
f(\sigma(i))=f(i)
$$
for every $i \in \{1, 2, 3, \dots, n\}$ form a subgroup of $S_{n}$,
we name this subgroup by $G_{f}$. In fact a $\sigma$ in $G_{f}$ is a
permutation among the inverse images of $f$ respectively, so $G_{f}$
isomorphs to the direct product of some symmetric groups.

For example, let $n=7$, $f(1)=f(2)=f(3)=f(4)=\alpha$,
$f(5)=f(6)=f(7)=\beta$, then $f^{-1}(\alpha)=\{1, 2, 3, 4\}$,
$f^{-1}(\beta)=\{5, 6, 7\}$, and $G_{f}\simeq S_{4}\times S_{3}$.

\begin{Definition}
For a given permutation $a \in S_{n}$, if the centralizer of $a$ in
$G_{f}$ is \{id\}, we call the permutation $a$ $f$-simple.
\end{Definition}

For example, $f(i)=i$ for every $i \in \{1, 2, 3, \dots, n\}$, then
$G_{f}=\{id\}$ and every permutation $a \in S_{n}$ is $f$-simple.
Another extreme case is when $f$ is a constant map, where
$G_{f}=S_{n}$ and then only $id \in S_{n}$ is $f$-simple.

It is clear that if a $\sigma \neq id \in G_{f}$ satisfies $\sigma
a=a \sigma$, then $f(a \sigma(i))=f(\sigma a(i))=f(a(i))$ for each
$i$. Furthermore, $f(a^{k} \sigma(i))=f(\sigma
a^{k}(i))=f(a^{k}(i))$ for each integer $k$.

For $\sigma \neq id$, there exists a number $i_{0}\in \{1, 2, 3,
\dots, n\}$ satisfying $\sigma(i_{0})\neq i_{0}$.
\begin{Lemma}
For a given permutation $a\in S_{n}$ and $\sigma\in G_{f}$ is a
centralizer element of $a$, suppose $\sigma(i_{0})\neq i_{0}$ for a
certain number $i_{0}\in \{1, 2, 3, \dots, n\}$.

If $\sigma(i_{0})=a^{d}(i_{0})$ for some $d>0$, then there exist a
integer $d_{0}>0$ and $k>1$, such that the cycle of permutation $a$
which including $i_{0}$ can be written as the form
$$
(i_{0}\cdots i_{d_{0}-1}i_{d_{0}}\cdots i_{kd_{0}-1})
$$
where $f(i_{cd_{0}+t})=f(i_{t})$ for $t=0,1,\dots ,d_{0}-1$ and
$c=0,1,\dots ,k-1$.
\end{Lemma}
\begin{pf}
Since $\sigma(i_{0})=a^{d}(i_{0})$ and $\sigma$ is a centralizer
element of $a$, we have
$$
\sigma a(i_{0})=a\sigma(i_{0})=a^{d+1}(i_{0})
$$
$$
\sigma a^{2}(i_{0})=a^{2}\sigma(i_{0})=a^{d+2}(i_{0})
$$
and so on. Following this and $\sigma\in G_{f}$ we obtain that
$$
f(a^{d}(i_{0}))=f(\sigma(i_{0}))=f(i_{0}),
$$
$$
f(a^{d+1}(i_{0}))=f(\sigma a(i_{0}))=f(a(i_{0}))
$$
$$
\cdots\cdots
$$
$$
f(a^{d+j}(i_{0}))=f(\sigma a^{j}(i_{0}))=f(a^{j}(i_{0}))
$$
$$
\cdots\cdots
$$
We then obtain that
$$
f(a^{hd+j}(i_{0}))=f(a^{j}(i_{0}))
$$
for every non-negative integers $h$ and $j$.

Assume $m=\textrm{min}\{c>0\mid a^{c}(i_{0})=i_{0}\}$,
$d_{0}=\textrm{gcd}(m,d)$. Thus $d_{0}<m$ since
$a^{d}(i_{0})=\sigma(i_{0})\neq i_{0}$. Then we can rewrite this
cycle of permutation $a$ to
$$
(i_{0}\cdots i_{d_{0}-1}i_{d_{0}}\cdots i_{kd_{0}-1}),
$$
where $k>1$ and $m=kd_{0}$.

Since $a^{m}(i_{0})=i_{0}$ we have
$$
f(a^{hm+j}(i_{0}))=f(a^{j}(i_{0}))
$$
for every non-negative integers $h$ and $j$. Furthermore, we can
write
$$
d_{0}=um+vd
$$
for some integers $u$ and $v$ since $d_{0}$ is the greatest common
divisor of $m$ and $d$. Following these we have
$$
f(a^{hd_{0}+j}(i_{0}))=f(a^{hum+hvd+j}(i_{0}))=f(a^{j}(i_{0}))
$$
for every non-negative integers $h$ and $j$. Following this we
finish our proof.
\end{pf}

\begin{Lemma}
For a given permutation $a\in S_{n}$ and $\sigma\in G_{f}$ is a
centralizer element of $a$, suppose $\sigma(i_{0})\neq i_{0}$ for a
certain number $i_{0}\in \{1, 2, 3, \dots, n\}$.

Denote the cycle of permutation $a$ which including $i_{0}$ as
$$
(i_{0}\cdots i_{m-1}).
$$
If $\sigma(i_{0})\neq a^{d}(i_{0})$ for every $d>0$, then the cycle
of permutation $a$ that including $\sigma(i_{0})$ has the form
$$
(j_{0}\cdots j_{m-1}),
$$
where $j_{0}=\sigma(i_{0})$. Therefore
$$
f(i_{k})=f(j_{k})
$$
for every $0\leq k\leq m-1$.
\end{Lemma}
\begin{pf}
Let $j_{k}=a^{k}\sigma(i_{0})$, $0\leq k\leq m-1$. Then these $m$
integers are all different with each other since $a\sigma=\sigma a$,
i.e. $(j_{0}\cdots j_{m-1})$ is a cycle of permutation $a$. Because
$\sigma(i_{0})\neq a^{d}(i_{0})$ for every $d>0$, it shows that the
cycle $(j_{0}\cdots j_{m-1})$ is different with $(i_{0}\cdots
i_{m-1})$.
\end{pf}

The inverse of \textbf{lemma 4.2.} and \textbf{lemma 4.3.} are also
right, that is
\begin{Lemma}
Suppose $a\in S_{n}$ is a permutation of $\{1, 2, 3, \dots, n\}$.

(i) If there is a cycle of the permutation $a$ has form
$$
(i_{0}\cdots i_{d-1}i_{d}\cdots i_{kd-1})
$$
where $d\geq 1, k>1$ and $f(i_{cd+t})=f(i_{t})$ for $t=0,1,\dots
,d-1$ and $c=0,1,\dots ,k-1$, then the permutation $a$ is not
$f$-simple.

(ii) If there are two cycles of the permutation $a$ have form
$$
(i_{0}\cdots i_{d-1})
$$
and
$$
(j_{0}\cdots j_{d-1})
$$
such that $d\geq1$ and
$$
f(i_{k})=f(j_{k})
$$
for every $0\leq k\leq d-1$, then the permutation $a$ is not
$f$-simple.
\end{Lemma}
\begin{pf}
For (i) the permutation
$$
\sigma=(i_{0}i_{d}i_{2d}\cdots
i_{(k-1)d})(i_{1}i_{d+1}i_{2d+1}\cdots i_{(k-1)d+1})\cdots
(i_{d-1}i_{2d-1}i_{3d-1}\cdots i_{kd-1})
$$
satisfies that $\sigma a=a \sigma$ and $\sigma\neq id$.

For (ii) the permutation
$$
\sigma=(i_{0}j_{0})(i_{1}j_{1})\cdots (i_{d-1}j_{d-1})
$$
satisfies that $\sigma a=a \sigma$ and $\sigma\neq id$.
\end{pf}

These lemmas give us an equivalent condition to the $f$-simple.
\begin{center}
\end{center}

Each cycle in a given permutation can be written as the form
$$
(i_{0}\cdots i_{d-1}i_{d}\cdots i_{kd-1})
$$
where $d\geq 1, k\geq 1$ and
$$
f(i_{cd+t})=f(i_{t})
$$
for $t=0,1,\dots ,d-1$ and $c=0,1,\dots ,k-1$. For example, let
$k=1$. We interest in how great $k$ can achieve for the given cycle.
If such a $k$ is the greatest one satisfying the condition above,
then the vectors
$$
(f(i_{0}),f(i_{1}),\dots, f(i_{d-1})),
$$
$$
(f(i_{1}),f(i_{2}),\dots, f(i_{d-1}), f(i_{0})),
$$
$$
\dots,
$$
$$
(f(i_{d-1}),f(i_{0}),\dots, f(i_{d-2}))
$$
are all different with each other, we call the set that formed by
these vectors the \textbf{$f$-kernel} of this given cycle. Note that
$f$-kernel is uniquely determined by the cycle and the map $f$.
\begin{center}
\end{center}

\begin{Proposition}
Suppose a set $G\subset S_{n}$ satisfying $\sigma G=G$ for all
$\sigma\in G_{f}$, then in $G$ the number of even no-$f$-simple
permutations is equal to the number of odd no-$f$-simple
permutations.
\end{Proposition}
\begin{pf}
Let $\tau\in G$, it can be written as the product of several
separated cycles uniquely. We classify these cycles by their
$f$-kernels and any class $C$ has the form below:
$$
C=(i_{10}\cdots i_{1(d-1)}i_{1d}\cdots i_{1(k_{1}d-1)})(i_{20}\cdots
i_{2(d-1)}i_{2d}\cdots i_{2(k_{2}d-1)})\cdots (i_{j0}\cdots
i_{j(d-1)}i_{jd}\cdots i_{j(k_{j}d-1)}),
$$
where $f(i_{s(cd+t)})=f(i_{1t})$ for $t=0,1,\dots ,d-1$, $s=1,\dots
,j$ and $c=0,1,\dots ,k_{s}-1$, i.e. these cycles in $C$ have same
$f$-kernel.

Let $\sigma_{C}$ be a permutation among $i_{s(cd+t)}$ which have
same value of $t$ respectively, then $\sigma_{C}\in G_{f}$. Since
$\sigma_{C}$ can be written as the product of several
transpositions, it is easily to prove that all cycles in $\sigma_{C}
C$ have the $f$-kernel same as of cycles in $C$.

Assume $\tau=C_{1}\cdot\cdots\cdot C_{m}$ where $C_{j}$ is the
product of the cycles of $\tau$ which have the same $f$-kernel,
$j=1, \cdots, m$. Since $\sigma G=G$ for all $\sigma\in G_{f}$, we
define the set $\{\prod_{j=1}^{m}\sigma_{C_{j}}C_{j}\}$ for any
$\sigma_{C_{j}}$ respect with $C_{j}$ to be the equivalent class of
$\tau$ over $G$. This equivalent relation over $G$ is well-defined
since all cycles in $\sigma_{C} C$ have the same $f$-kernel as of
cycles in $C$.

Note that all $\sigma_{C}$ form a subgroup of $G_{f}$ which
isomorphs to $(S_{\sum_{s=1}^{j}k_{s}})^{d}$ with
$S_{\sum_{s=1}^{j}k_{s}}$ the symmetric group of degree
$\sum_{s=1}^{j}k_{s}$, then following $G_{f}G=G$ the equivalent
class of $\tau$ over $G$ is $\tau$ left multiplied by such subgroup
of $G_{f}$ which isomorphs to the product of
$(S_{\sum_{s=1}^{j}k_{s}})^{d}$. If $\tau$ is no-$f$-simple
permutation, then some $\sum_{s=1}^{j}k_{s}>1$ by \textbf{lemma
4.2.} and \textbf{lemma 4.3.}. Thus in the equivalent class of
$\tau$ over $G$, the number of even no-$f$-simple permutations is
equal to the number of odd no-$f$-simple permutations. This is what
we want to prove.
\end{pf}

Assume $a$ is an $f$-simple permutation, then for any $\sigma\in
G_{f}$, the conjugation $\sigma a\sigma^{-1}$ is not equal to $a$
when $\sigma\neq id$. Therefore $|\{\sigma a\sigma^{-1}\mid\sigma\in
G_{f}\}|=|G_{f}|$, and following \textbf{proposition 4.5.} we have

\begin{Proposition}
Suppose a set $G\subset S_{n}$ satisfying $\sigma G=G$ and $\sigma a
\sigma^{-1}\in G$ for all $\sigma\in G_{f}$ and all $f$-simple
permutation $a\in G$. If the number of even permutations in $G$ is
equal to the number of odd permutations in $G$, then in $G$ the
number of such conjugate classes
$$
\{\sigma a\sigma^{-1}\mid\sigma\in G_{f}\}
$$
where $a\in G$ is even $f$-simple permutation, is equal to the
number of such conjugate classes where $a\in G$ is odd $f$-simple
permutation.
\end{Proposition}
\begin{center}
\end{center}

\section{\bf Permutation of $w$}

By the reason in \textbf{section 3}, a term as (7) determines such a
table in \textbf{section 3}. Recalling the definition of
$\overline{w_{j}}$, a term as (7) also determines a set
$W=\{\overline{w_{1}},\cdots , \overline{w_{s}}\}$ of $m$-dimension
vectors(same elements do not combine). We say that two terms as (7)
are equivalent, if they determine same set of $W$. Obviously this is
an equivalent relationship, and for each term in the same equivalent
class, following equation (8) we will see that they have the same
part of
$\prod_{i=1}^{m}(\chi(a_{i})^{\sum_{j}\sigma((p^{s_{j}}-1)r_{ji})}\pi^{\sum_{j}\sigma((p^{s_{j}}-1)r_{ji})})$
since $\sum_{j}\sigma((p^{s_{j}}-1)r_{ji})$ is equal to the sum of
all $k_{ji}[t]$ that in the table.

Besides, following Stickelberger theorem the part
$$
\prod_{j}(-\prod_{i=1}^{m}\prod_{t=0}^{s_{j}-1}\Gamma_{p}(\{p^{t}r_{ji}\}))
$$
$\ \textrm{mod}\  p$ is equal to
$\prod_{j}\prod_{i=1}^{m}\prod_{t=0}^{s_{j}-1}\frac{1}{k_{ji}[t]!}$
or its inverse, it depends on whether the number of index $j$ is
even or not, i.e. it depends on whether the number of blocks in the
corresponding table is even or not. This means all terms in the same
equivalent class are at most different with a sign after $\
\textrm{mod}\ p$.

For a fixed set $W=\{\overline{w_{1}},\dots , \overline{w_{s}}\}$
corresponding to some term in $T(c_{s})$, let $f$ be the injective
map from $\{w_{1},\dots , w_{s}\}$ to $W$ defined by
$$
f:w_{t}\rightarrow \overline{w_{t}},
$$
and consider $S_{s}$ the symmetric group of $\{w_{1},\dots ,
w_{s}\}$. Since a permutation $a\in S_{s}$ can be written as the
product of several separated cycles uniquely, when replace them by
their $f$-values the permutation $a$ determines a table as the form
mentioned in \textbf{section 3} directly, and each of cycles
transforms to a block. Following \textbf{lemma 4.2.}, \textbf{lemma
4.3.} and \textbf{lemma 4.4.} we can easily see that a table is
corresponding to a term as (7) if and only if it determines an
$f$-simple permutation of $S_{s}$ and satisfies the relation (9) in
\textbf{proposition 3.4.}. Define
$$
G=\{a\in S_{s}\mid \textrm{the table corresponding to } a \textrm{
satisfying (9) }\}.
$$

Suppose there are $k$ distinct $u$-values of $\overline{w_{t}}\in
W$, denote by $u_{1}, \cdots, u_{k}$. Let
$$
W_{i}=\{w_{t} \mid \overline{w_{t}}\in W \textrm{ and } \textrm{the
} u-\textrm{value } \textrm{of }\ \overline{w_{t}}\ \textrm{ is
equal to } u_{i}\}
$$
and
$$
s_{i}=\mid W_{i}\mid
$$
for $i=1, \cdots, k$. Then
\begin{Lemma}
Let $S_{s_{i}}$ be the symmetric group of $W_{i}$, then for any
$a\in G$, $a^{-1}G=\prod_{i=1}^{k}S_{s_{i}}$.
\end{Lemma}
\begin{pf}
For any $b\in G$ and any $w\in\{w_{1},\dots , w_{s}\}$, by
definition of $G$, $\overline{w}$'s $u$-value is equal to
$\overline{b(w)}$'s $v$-value. Thus $\overline{a^{-1}b(w)}$'s
$u$-value is equal to $\overline{w}$'s $u$-value. This shows that
$a^{-1}G\subset\prod_{i=1}^{k}S_{s_{i}}$.

Besides, for any $c\in \prod_{i=1}^{k}S_{s_{i}}$ and any
$w\in\{w_{1},\dots , w_{s}\}$, $\overline{c(w)}$'s $u$-value is
equal to $\overline{w}$'s $u$-value. Thus $ac\in G$. This shows that
$\prod_{i=1}^{k}S_{s_{i}}\subset a^{-1}G$.

Then we complete the proof.
\end{pf}

If there are two elements of $W$ having same $u$-value, then some of
$s_{i}>1$ and by \textbf{lemma 5.1.} the number of even permutations
is equal to the number of odd permutations in $G$.

Since any permutation can be written as a product of several
transpositions, it is easily to check that $G_{f}G=GG_{f}=G$, thus
it satisfies the condition of \textbf{proposition 4.5.}.

For any $\sigma\in G_{f}$ and $a\in G$ an $f$-simple permutation,
$\sigma a \sigma^{-1}$ and $a$ are corresponding to the same term as
(7), and so its inverse, i.e. for any two $f$-simple permutations
$a, b\in G$, if they are corresponding to the same term as (7) by
replacing $f$-values, then there exist a $\sigma\in G_{f}$ such that
$b=\sigma a \sigma^{-1}$. Following \textbf{proposition 4.6.} we
obtain that there is a half number of terms in (8) for which the
part
$\prod_{j}(-\prod_{i=1}^{m}\prod_{t=0}^{s_{j}-1}\Gamma_{p}(\{p^{t}r_{ji}\}))$
$\ \textrm{mod}\  p$ is equal to
$\prod_{j}\prod_{i=1}^{m}\prod_{t=0}^{s_{j}-1}\frac{1}{k_{ji}[t]!}$
and the same number terms for which the part
$\prod_{j}(-\prod_{i=1}^{m}\prod_{t=0}^{s_{j}-1}\Gamma_{p}(\{p^{t}r_{ji}\}))$
$\ \textrm{mod}\  p$ is equal to
$-\prod_{j}\prod_{i=1}^{m}\prod_{t=0}^{s_{j}-1}\frac{1}{k_{ji}[t]!}$.
By the reason above, we obtain that
\begin{Proposition}
Assume that $p\geq \sum_{i=1}^{m}d_{i}$. If there are two vectors
$\overline{w}$ of $W$ corresponding to the same value of $u$, then
the sum of all terms as (7) in the same equivalent class has the
$\textrm{ord}_{p}$-value not smaller than
$1+\frac{1}{p-1}\sum_{j}\sum_{i=1}^{m}\sigma((p^{s_{j}}-1)r_{ji})$.
\end{Proposition}

\begin{Corollary}
In \textbf{proposition 5.2.}, the lower bound can be instead by
$1+\frac{1}{d}\sum_{j}\sum_{t=0}^{s_{j}-1}u_{j}[t]$.
\end{Corollary}

\begin{pf}
By definition, $\sum_{i=1}^{m}d_{i}k_{ji}[t]=u_{j}[t]p-v_{j}[t]$,
thus  $\sum_{i=1}^{m}k_{ji}[t]\geq
\frac{1}{d}\sum_{i=1}^{m}d_{i}k_{ji}[t]=\frac{1}{d}(u_{j}[t]p-v_{j}[t])$.
So
\begin{eqnarray*}
\frac{1}{p-1}\sum_{j}\sum_{i=1}^{m}\sigma((p^{s_{j}}-1)r_{ji})
&=&\frac{1}{p-1}\sum_{j}\sum_{i=1}^{m}\sigma(k_{ji}) \\
&=&\frac{1}{p-1}\sum_{j}\sum_{i=1}^{m}\sum_{t=0}^{s_{j}-1}k_{ji}[t]\\
&=&\frac{1}{p-1}\sum_{j}\sum_{t=0}^{s_{j}-1}\sum_{i=1}^{m}k_{ji}[t]\\
&\geq&
\frac{1}{p-1}\sum_{j}\sum_{t=0}^{s_{j}-1}\frac{1}{d}(u_{j}[t]p-v_{j}[t])\\
&=&\frac{1}{d}\sum_{j}\sum_{t=0}^{s_{j}-1}u_{j}[t],
\end{eqnarray*}
the last equation follows from the relationship (9). Then we have
proved the corollary.
\end{pf}
\begin{center}
\end{center}

If $u_{j}[t]=0$ for some $j$ and $t$, then $v_{j}[t]$ must also be
the value $0$. By \textbf{proposition 3.4.} we can obtain that the
values $u$, $v$ in every column vectors of this block must be $0$,
it means that this block has only one column. Because any table of
term has at most one such a block, thus
$\sum_{j}\sum_{t=0}^{s_{j}-1}u_{j}[t]\geq s-1$, that is

\begin{Corollary}
In \textbf{proposition 5.2.}, the lower bound can be instead by
$1+\frac{s-1}{d}$.
\end{Corollary}

Following \textbf{corollary 5.4.} we have proved the first part of
\textbf{proposition 3.5.}.

By the proof of \textbf{corollary 5.3.}, for any $c\in T(c_{s})$,
$$
\textrm{ord}_{p}c\geq
\frac{1}{d}\sum_{j}\sum_{t=0}^{s_{j}-1}u_{j}[t].
$$
If any two of $u_{j}[t]$ are all different, then
$$
\textrm{ord}_{p}c\geq \frac{1}{d}\sum_{u=0}^{s-1}u=
\frac{s(s-1)}{2d}.
$$

Besides, the Newton polygon of $L(f,T)$ is symmetric in the sense
that for every slope segment $\alpha$ there is a slope segment
$1-\alpha$ of the same horizontal length, we should only determine
half number slope segments of the Newton polygon of $L(f,T)$. In
other words, for $L^{*}(f,T)= (1-T)L(f,T)$, we should only consider
those coefficients $c_{s}$ where $s$ satisfies
$s-1\leq\frac{d-1}{2}$.

Thus, to prove the last part of \textbf{proposition 3.5.}, we should
only show that
$$
\frac{s(s-1)}{2d}< 1+\frac{s-1}{d},
$$
which is equivalent to $(s-1)(s-2)<2d$.
\begin{center}
\end{center}

\section{\bf General calculation of $\textrm{ord}_{p}c_{s}$ for $(s-1)(s-2)<2d$}

Let $m>1$, $f(x)= \sum_{i=1}^{m}a_{i}x^{d_{i}}$ be a polynomial,
where $0<d_{1}<\cdots<d_{m}= d$ and $a_{i}\in \textbf{F}_{p}^{*}$,
$a_{m}=1$. Assume $p\geq \sum_{i=1}^{m}d_{i}$. Since such a linear
transformation $ax+b$ $(a,b\in \textbf{F}_{q}, a\neq 0 \
\textrm{mod}\  p)$ of $x$ does not change the $L$-function, we can
also assume $d_{m-1}<d-1$.

By \textbf{proposition 3.5.} we should consider such $c\in T(c_{s})$
that corresponds to $s$ different $u$-values beginning  from $0, 1,
\cdots, s-1$ respectively. Besides, for \textbf{proposition 3.4.}
the $v$-values of the table corresponding to $c$ are equal to its
$u$-values respectively. To determine a $c\in T(c_{s})$ is
equivalent to determine the relevant table. Therefore we should
consider these $s$ equations:
\begin{equation}
\sum_{i=1}^{m}d_{i}k_{ji}[t]=u_{j}[t]p-v_{j}[t]
\end{equation}
for all $j$ and $t$ such that all $k_{ji}[t]$ are non-negative
integers.

Recall that
$$
\textrm{ord}_{p}c=\frac{1}{p-1}\sum_{j}\sum_{i=1}^{m}\sum_{t=0}^{s_{j}-1}k_{ji}[t],
$$
for a series given positive integers $r_{j}[t]$, we insert $s$
equations
\begin{equation}
\sum_{i=1}^{m}k_{ji}[t]=r_{j}[t]
\end{equation}
into (14). Thus, to calculate $\textrm{ord}_{p}c_{s}$ we should
calculate $\textrm{ord}_{p}(\sum_{\textrm{ord}_{p}c=\frac{r}{p-1}}
c)$ satisfying (14) and (15) in order of
$r=\sum_{j}\sum_{t=0}^{s_{j}-1}r_{j}[t]$ from small to large.

If some $k_{ji}[t]\geq d$ with $i<m$, then we can use
$k_{jm}[t]+d_{i}$ and $k_{ji}[t]-d$ instead of $k_{jm}[t]$ and
$k_{ji}[t]$ satisfying (14) and small
$r_{j}[t]=\sum_{i=1}^{m}k_{ji}[t]$ into
$\sum_{i=1}^{m}k_{ji}[t]-d+d_{i}$. Therefore, to make $r$ smallest
all $k_{ji}[t]$ with $i<m$ should be smaller than $d$.

Combine (14) and (15) we get
\begin{equation}
\sum_{i=1}^{m-1}(d-d_{i})k_{ji}[t]=dr_{j}[t]-u_{j}[t]p+v_{j}[t].
\end{equation}

Define $C(r; u, v)$ the set of all non-negative integral solutions
$[h_{1}, h_{2}, \cdots, h_{d-2}]$ of
$$
\sum_{i=1}^{m-1}(d-d_{i})k_{i}=dr-up+v
$$
with $h_{d_{i}}=k_{i}$ for $i=1, \cdots, m-1$ and $h_{j}=0$ for all
indexes $j\neq d_{1}, \cdots, d_{m}-1$.

Suppose $u_{1}, \cdots, u_{s}$ are $s$ distinct non-negative
integers, and $\sigma\in S_{s}$ is a permutation on $\{1, \cdots,
s\}$. Let $r=\sum_{i=1}^{s}r_{i}$. When we select $s$ solutions in
$C(r_{1}; u_{1}, u_{\sigma(1)}), \cdots,$ $C(r_{s}; u_{s},
u_{\sigma(s)})$ respectively, by \textbf{proposition 3.4.} they
construct a table uniquely, which is corresponding to a term $c\in
T(c_{s})$ with $\textrm{ord}_{p}c=\frac{r}{p-1}$.

Recall \textbf{section 5}, the part
$$
\prod_{j}(-\prod_{i=1}^{m}\prod_{t=0}^{s_{j}-1}\Gamma_{p}(\{p^{t}r_{ji}\}))
$$
$\textrm{mod}\  p$ is equal to
$\prod_{j}\prod_{i=1}^{m}\prod_{t=0}^{s_{j}-1}\frac{1}{k_{ji}[t]!}$
or its inverse, it depends on whether the number of blocks in the
corresponding table is even or not, i.e. it depends on whether the
permutation
$$
\left(
\begin{array}{cccc}
u_{1} & u_{2} & \cdots & u_{s} \\
v_{1} & v_{2} & \cdots & v_{s} \\
\end{array}
\right)
$$
is even or odd. We can also calculate
$\prod_{j}\prod_{i=1}^{m}\prod_{t=0}^{s_{j}-1}\frac{1}{k_{ji}[t]!}$
column by column in the table, i.e. a column $[k_{1}, \cdots,
k_{m}]^{T}=[h_{d_{1}}, \cdots, h_{d_{m-1}},
r-\sum_{i=1}^{m-1}h_{d_{i}}]^{T}$ with $[h_{1}, \cdots, h_{d-2}]\in
C(r; u, v)$ supplies the part
$$
(\prod_{i=1}^{m-1}\frac{1}{k_{i}!})\cdot
\frac{1}{(r-\sum_{i=1}^{m-1}k_{i})!} \ \textrm{mod}\  p.
$$

Define
$$
F(r; u, v)=\sum_{[h_{1}, \cdots, h_{d-2}]\in C(r; u,
v)}(\prod_{i=1}^{m-1}\frac{1}{k_{i}!}\chi(a_{i})^{k_{i}})\cdot
\frac{1}{(r-\sum_{i=1}^{m-1}k_{i})!}
$$
and
$$
F_{r}^{s}=\sum_{\sigma\in
S_{s}}\sum_{\sum_{i=1}^{s}r_{i}=r}\sum_{u_{1}, \cdots, u_{s}
\textrm{ are all
distinct}}\textrm{sign}(\sigma)\prod_{i=1}^{s}F(r_{i}; u_{i},
u_{\sigma(i)}).
$$
If we define $O(c_{s})=\{c\in T(c_{s})\mid \textrm{all the }
u-\textrm{values corresponding to } c\ \textrm{are distinct each
other}. \}$, then we have
$$
\pi^{-r}(\sum_{c\in O(c_{s}),\
\textrm{ord}_{p}c=\frac{r}{p-1}}c)\equiv F_{r}^{s} \ \textrm{mod}\
p.
$$

Note that when $r$ is small, the $\sum_{u_{1}, \cdots, u_{s}
\textrm{ are all distinct}}$ of $F_{r}^{s}$ only contains the case
$\{u_{1}, \cdots, u_{s}\}=\{0, \cdots, s-1\}$.

By \textbf{proposition 3.5.}, to determine $\textrm{ord}_{p}c_{s}$,
we should only begin to calculate $F_{r}^{s}$ in order of $r$ from
small to large until for the first $r$ satisfying $F_{r}^{s}\neq 0\
\textrm{mod}\  p$. we also see that all $r_{i}$ are disjointed, then
to make $r$ smallest we should only make each $r_{i}$ smallest.

By \textbf{proposition 3.5.} again we have

\begin{Theorem}
Suppose $F_{r}^{s}=0$ when $r<R$ and $F_{R}^{s}\neq 0$, and denote
$\lambda_{s}=\frac{R}{p-1}$. If $\lambda_{s}<1+\frac{s-1}{d}$ then
$$
\textrm{ord}_{p}c_{s}=\lambda_{s}.
$$
\end{Theorem}

\begin{Remark}
This theorem give a feasible method to calculate Newton polygon of
$L$-function in sense of Sperber's theorem on case $d=3$ (see [16]).
\end{Remark}

\begin{center}
\end{center}

\section{\bf General calculation for $d=3, 4, 6$}

Let $\{\omega_{0},\dots,\omega_{d-1}\}$ be the set of reciprocal
roots of $L^{*}(f,T)$ satisfying
$$
0=\textrm{ord}_{p}\omega_{0}\leq\dots\leq\textrm{ord}_{p}\omega_{d-1}
$$
then $\{\omega_{1},\dots,\omega_{d-1}\}$ is the set of reciprocal
roots of $L(f,T)$.

If $p\equiv 1\ \textrm{mod}\  d$, then $[0, 0, \cdots, 0]\in
C(\frac{u(p-1)}{d}; u, u)$ ($u=0, 1, \cdots, s-1$) is the unique
solution satisfying (16) that achieves the lower bound
$\sum_{u=0}^{s-1}\frac{u(p-1)}{d}$ of $r=\sum_{i=1}^{s}r$. Thus
$$
\textrm{ord}_{p}\omega_{s}=\frac{s-1}{d}
$$
for all $s$. This is the $1$-dimension case of A.S. conjecture.

Generally, by the discussion above, we first consider the terms in
$c_{s}$ of which $W$-set is corresponding to the $u$-values set as
$\{0, 1, \dots, s-1\}$ and calculate the sum of terms among them
which has the smallest $\textrm{ord}_{p}$-value. if the sum not
increase the $\textrm{ord}_{p}$-value, then it is the
$\textrm{ord}_{p}$-value of $c_{s}$; otherwise the sum of their
$\Gamma_{p}$-part is equal to $0\ \textrm{mod}\  p$, therefore the
$\textrm{ord}_{p}$-value at least increase $1$ and then we calculate
the sum of terms among them having the next bigger
$\textrm{ord}_{p}$-value and so on. Note that
$\textrm{ord}_{p}c_{1}$ is always equal to 0, we calculate
$\textrm{ord}_{p}c_{2}$ first of all, it is corresponding to the
first slope segment of Newton polygon of $L(f,T)$.

Since the Newton polygon of $L(f,T)$ is symmetric in the sense that
for every slope segment $\alpha$ there is a slope segment $1-\alpha$
of the same horizontal length, $x-1=2y$ lies upon the Newton
polygon. Besides, the points $(s, 1+\frac{s-1}{d})$ are all on the
line $y=1+\frac{x-1}{d}$. Thus if $(1+\frac{d-1}{2},
1+\frac{d-1}{2d})$ is over line $x-1=2y$, i.e. if $d\leq 6$, then
following \textbf{proposition 3.5.} we have determined the Newton
polygon by \textbf{theorem 6.1.}.

\begin{Theorem}
Let $f(x)=x^{3}+a_{1}x$, $a_{1}\neq0$ and $p>3$, $p\equiv 2\
\textrm{mod}\  3$. Thus
$$
\textrm{ord}_{p}\omega_{1}=\frac{p+1}{3(p-1)},
\textrm{ord}_{p}\omega_{2}=1-\frac{p+1}{3(p-1)}.
$$
\end{Theorem}

\begin{pf}
Let $s=2$, then we should consider (16) when $u=v=1$. Since
$\frac{p+1}{3}$ is an integer, thus $r=\frac{p+1}{3}$ is the lower
bound of $r$ while $C(\frac{p+1}{3}; 1, 1)=\{[1]\}$. So
$F(\frac{p+1}{3}; 1, 1)=\frac{1}{(\frac{p-2}{3})!}\chi(a_{1})$ and
$F_{\frac{p+1}{3}}^{2}=\frac{1}{(\frac{p-2}{3})!}\chi(a_{1})\neq 0\
\textrm{mod}\ p$. Thus
$\textrm{ord}_{p}c_{2}=\lambda_{2}=\frac{p+1}{3(p-1)}$ and we finish
the proof.
\end{pf}

\begin{Remark}
Assume $f(x)=x^{3}+a_{2}x^{2}+a_{1}x$, $a_{1}\neq0, a_{2}\neq0$,
$p>3$, $p\equiv 2\ \textrm{mod}\  3$. Then the lower bound of $r$ is
$\frac{p+1}{3}$ and $C(\frac{p+1}{3}; 1, 1)=\{[1, 0], [0, 2]\}$. So
\begin{eqnarray*}
F_{\frac{p+1}{3}}^{2} &=&
\frac{1}{(\frac{p-2}{3})!}\chi(a_{1})+\frac{1}{2!}\cdot\frac{1}{(\frac{p-2}{3}-1)!}\chi(a_{2})^{2}
\\
&\equiv& \frac{1}{3(\frac{p-2}{3})!}(3\chi(a_{1})-\chi(a_{2})^{2})\
\textrm{mod}\ p.
\end{eqnarray*}
If $F_{\frac{p+1}{3}}^{2}\equiv 0\ \textrm{mod}\ p$, i.e.
$3\chi(a_{1})-\chi(a_{2})^{2}=0\ \textrm{mod}\ p$, we will consider
the second smallest value of $r$, probably is $\frac{p+1}{3}+1$. For
$C(\frac{p+1}{3}+1; 1, 1)=\{[2, 1], [1, 3], [0, 5]\}$,
$$
F_{\frac{p+1}{3}+1}^{2}=\chi(a_{1})^{2}\chi(a_{2})\cdot\frac{1}{2!}\cdot\frac{1}{(\frac{p-2}{3}-1)!}+\chi(a_{1})\chi(a_{2})^{3}\cdot\frac{1}{3!}\cdot\frac{1}{(\frac{p-2}{3}-2)!}+\chi(a_{2})^{5}\cdot\frac{1}{5!}\cdot\frac{1}{(\frac{p-2}{3}-3)!}.
$$
Since $3\chi(a_{1})-\chi(a_{2})^{2}=0\ \textrm{mod}\ p$, we have
\begin{eqnarray*}
F_{\frac{p+1}{3}+1}^{2}&\equiv& \chi(a_{2})^{5}\cdot
(\frac{1}{3^{2}}\cdot\frac{1}{2!}\cdot\frac{1}{(\frac{p-2}{3}-1)!}+\frac{1}{3}\cdot\frac{1}{3!}\cdot\frac{1}{(\frac{p-2}{3}-2)!}+\frac{1}{5!}\cdot\frac{1}{(\frac{p-2}{3}-3)!}) \\
&\equiv& \frac{\chi(a_{2})^{5}}{(\frac{p-2}{3})!}\cdot
(\frac{1}{3^{2}}\cdot\frac{1}{2!}\cdot\frac{-2}{3}+\frac{1}{3}\cdot\frac{1}{3!}\cdot(\frac{-2}{3})(\frac{-2}{3}-1)+\frac{1}{5!}\cdot(\frac{-2}{3})(\frac{-2}{3}-1)(\frac{-2}{3}-2)) \\
&\equiv& 0\ \textrm{mod}\ p,
\end{eqnarray*}
and then we will consider the third smallest value of $r$, and so
on. This is similar as [16].

However, if we changing $x$ by $x-\frac{1}{3}a_{2}$, following
$3\chi(a_{1})-\chi(a_{2})^{2}= 0 \ \textrm{mod}\  p$ we have
$$
f(x-\frac{1}{3}a_{2})=x^{3}-\frac{1}{27}a_{2}^{3},
$$
it is diagonal case, so that
$$
\textrm{ord}_{p}\omega_{1}=\textrm{ord}_{p}\omega_{2}=\frac{1}{2}.
$$
It is better than [16].
\end{Remark}

\begin{Theorem}
Let $f(x)=x^{4}+a_{2}x^{2}+a_{1}x$ and $p>6$, $p\equiv 3\
\textrm{mod}\  4$.

If $a_{2}\neq 0$, then
$$
\textrm{ord}_{p}\omega_{1}=\frac{p+1}{4(p-1)},
\textrm{ord}_{p}\omega_{2}=\frac{1}{2},
\textrm{ord}_{p}\omega_{3}=1-\frac{p+1}{4(p-1)};
$$

If $a_{2}=0$ and $a_{1}\neq 0$, then
$$
\textrm{ord}_{p}\omega_{1}=\frac{p+5}{4(p-1)},
\textrm{ord}_{p}\omega_{2}=\frac{1}{2},
\textrm{ord}_{p}\omega_{3}=1-\frac{p+5}{4(p-1)}.
$$
\end{Theorem}

\begin{pf}
Let $s=2$, then we should consider (16) when $u=v=1$.

Suppose $a_{2}\neq 0$, then $r=\frac{p+1}{4}$ is the lower bound of
$r$ while $C(\frac{p+1}{4}; 1, 1)=\{[0, 1]\}$. So $F(\frac{p+1}{4};
1, 1)=\frac{1}{(\frac{p-3}{4})!}\chi(a_{2})$ and
$F_{\frac{p+1}{4}}^{2}=\frac{1}{(\frac{p-3}{4})!}\chi(a_{2})\neq 0 \
\textrm{mod}\ p$. Thus
$\textrm{ord}_{p}c_{2}=\lambda_{2}=\frac{p+1}{4(p-1)}$. We finish
the proof of the first case.

Suppose $a_{2}=0, a_{1}\neq 0$, then $r=\frac{p+5}{4}$ is the lower
bound of $r$ while $C(\frac{p+5}{4}; 1, 1)=\{[2,0]\}$. So
$F(\frac{p+5}{4}; 1,
1)=\frac{1}{2!}\cdot\frac{1}{(\frac{p-3}{4})!}\chi(a_{1})^{2}$ and
$F_{\frac{p+5}{4}}^{2}=\frac{1}{2!}\cdot\frac{1}{(\frac{p-3}{4})!}\chi(a_{1})^{2}\neq
0 \ \textrm{mod}\ p$. Thus
$\textrm{ord}_{p}c_{2}=\lambda_{2}=\frac{p+5}{4(p-1)}$. We finish
the proof of the last case.
\end{pf}

\begin{Remark}
In fact the first part in this theorem has two cases: $a_{1}=0$ and
$a_{1}\neq 0$. Since $C(\frac{p+1}{4}; 1, 1)=\{[0, 1]\}$ does not
depend on whether $a_{1}=0$, it does not affect the calculation of
$F_{\frac{p+5}{4}}^{2}$. That is why we defined $C(r; u, v)$ like
that.
\end{Remark}
\begin{center}
\end{center}

To simply describe the results, we use $a_{i}$ directly to instead
$\chi(a_{i})$.

\begin{Theorem}
Let $f(x)=x^{6}+a_{4}x^{4}+a_{3}x^{3}+a_{2}x^{2}+a_{1}x$ and $p\geq
6+\sum_{a_{i}\neq 0}i$, $p\equiv-1\ \textrm{mod}\  6$.

\textrm{(i)} If $a_{4}\neq 0$, $3a_{3}^{2}+a_{4}^{3}-3a_{2}a_{4}\neq
0$, then
\begin{eqnarray*}
&\ &\textrm{ord}_{p}\omega_{1}=\frac{p+1}{6(p-1)},
\textrm{ord}_{p}\omega_{2}=\frac{p+1}{3(p-1)},
\textrm{ord}_{p}\omega_{3}=\frac{1}{2},
\textrm{ord}_{p}\omega_{4}=1-\frac{p+1}{3(p-1)},\\
&\ &\textrm{ord}_{p}\omega_{5}=1-\frac{p+1}{6(p-1)}.
\end{eqnarray*}

\textrm{(ii)}If $p\geq 17$ (resp. $p=11$), $a_{4}\neq 0$,
$3a_{3}^{2}+a_{4}^{3}-3a_{2}a_{4}=0$,
$7a_{3}a_{4}^{3}+12a_{3}^{3}-12a_{1}a_{4}^{2}\neq 0$ (resp.
$-4a_{1}^{2}a_{4}^{4}+2a_{3}^{2}a_{4}^{6}-a_{3}^{4}a_{4}^{3}-4a_{3}^{6}+a_{1}a_{3}a_{4}^{5}-3a_{1}a_{3}^{3}a_{4}^{2}\neq
0$), then
\begin{eqnarray*}
&\ &\textrm{ord}_{p}\omega_{1}=\frac{p+1}{6(p-1)},
\textrm{ord}_{p}\omega_{2}=\frac{p+4}{3(p-1)},
\textrm{ord}_{p}\omega_{3}=\frac{1}{2},
\textrm{ord}_{p}\omega_{4}=1-\frac{p+4}{3(p-1)}, \\
&\ &\textrm{ord}_{p}\omega_{5}=1-\frac{p+1}{6(p-1)}.
\end{eqnarray*}

\textrm{(iii)} If $p\geq 29$(resp. $p=11$, $p=17$, $p=23$),
$a_{4}\neq 0$, $3a_{3}^{2}+a_{4}^{3}-3a_{2}a_{4}=0$,
$7a_{3}a_{4}^{3}+12a_{3}^{3}-12a_{1}a_{4}^{2}=0$(resp.
$-4a_{1}^{2}a_{4}^{4}+2a_{3}^{2}a_{4}^{6}-a_{3}^{4}a_{4}^{3}-4a_{3}^{6}+a_{1}a_{3}a_{4}^{5}-3a_{1}a_{3}^{3}a_{4}^{2}=0$,
$7a_{3}a_{4}^{3}+12a_{3}^{3}-12a_{1}a_{4}^{2}=0$,
$7a_{3}a_{4}^{3}+12a_{3}^{3}-12a_{1}a_{4}^{2}=0$),
$48384a_{3}^{8}+225a_{4}^{12}+5600a_{4}^{9}a_{3}^{2}+41888a_{4}^{6}a_{3}^{4}+80640a_{4}^{3}a_{3}^{6}\neq
0$(resp.
$5+5a_{3}^{4}a_{4}^{4}+2a_{3}^{6}a_{4}-2a_{3}^{8}a_{4}^{8}-2a_{4}^{5}+5a_{3}^{2}a_{4}^{2}-2a_{3}^{4}a_{4}^{9}-2a_{3}^{6}a_{4}^{6}+2a_{3}^{8}a_{4}^{3}\neq
0$,
$-4a_{3}^{8}+8a_{4}^{5}a_{3}^{4}+7a_{3}^{5}a_{4}^{4}-2a_{3}^{7}a_{4}^{2}\neq
0$,
$-10a_{3}^{8}+11a_{4}^{12}+8a_{3}^{2}a_{4}^{9}-11a_{3}^{4}a_{4}^{6}-9a_{3}^{6}a_{4}^{3}\neq
0$), then
\begin{eqnarray*}
&\ &\textrm{ord}_{p}\omega_{1}=\frac{p+1}{6(p-1)},
\textrm{ord}_{p}\omega_{2}=\frac{p+7}{3(p-1)}(when \ p=11 \ here \
is \
\frac{1}{2} \ instead),\\
&\ &\textrm{ord}_{p}\omega_{3}=\frac{1}{2},
\textrm{ord}_{p}\omega_{4}=1-\frac{p+7}{3(p-1)}(when \ p=11 \ here \
is
\ \frac{1}{2} \ instead),\\
&\ &\textrm{ord}_{p}\omega_{5}=1-\frac{p+1}{6(p-1)}.
\end{eqnarray*}

\textrm{(iv)} If $p\geq 29$(resp. $p=17$, $p=23$), $a_{4}\neq 0$,
$3a_{3}^{2}+a_{4}^{3}-3a_{2}a_{4}=0$,
$7a_{3}a_{4}^{3}+12a_{3}^{3}-12a_{1}a_{4}^{2}=0$(resp.
$7a_{3}a_{4}^{3}+12a_{3}^{3}-12a_{1}a_{4}^{2}=0$,
$7a_{3}a_{4}^{3}+12a_{3}^{3}-12a_{1}a_{4}^{2}=0$),
$48384a_{3}^{8}+225a_{4}^{12}+5600a_{4}^{9}a_{3}^{2}+41888a_{4}^{6}a_{3}^{4}+80640a_{4}^{3}a_{3}^{6}=
0$(resp.
$-4a_{3}^{8}+8a_{4}^{5}a_{3}^{4}+7a_{3}^{5}a_{4}^{4}-2a_{3}^{7}a_{4}^{2}=
0$,
$-10a_{3}^{8}+11a_{4}^{12}+8a_{3}^{2}a_{4}^{9}-11a_{3}^{4}a_{4}^{6}-9a_{3}^{6}a_{4}^{3}=
0$), then
\begin{eqnarray*}
&\ &\textrm{ord}_{p}\omega_{1}=\frac{p+1}{6(p-1)},
\textrm{ord}_{p}\omega_{2}=\frac{p+10}{3(p-1)}(when \ p=17 \ here \
is \
\frac{1}{2} \ instead),\\
&\ &\textrm{ord}_{p}\omega_{3}=\frac{1}{2},
\textrm{ord}_{p}\omega_{4}=1-\frac{p+10}{3(p-1)}(when \ p=17 \ here
\ is
\ \frac{1}{2} \ instead),\\
&\ &\textrm{ord}_{p}\omega_{5}=1-\frac{p+1}{6(p-1)}.
\end{eqnarray*}

\textrm{(v)} If $a_{3}=a_{1}=0$, $3a_{2}=a_{4}^{2}\neq 0$, then
$$
\textrm{ord}_{p}\omega_{1}=\frac{p+1}{6(p-1)},
\textrm{ord}_{p}\omega_{2}=\textrm{ord}_{p}\omega_{3}=\textrm{ord}_{p}\omega_{4}=\frac{1}{2},
\textrm{ord}_{p}\omega_{5}=1-\frac{p+1}{6(p-1)}.
$$

\textrm{(vi)} If $a_{4}=0$, $a_{3}\neq 0$,
$a_{2}^{2}+2a_{1}a_{3}\neq 0$, and $a_{1}\neq 0$ or $a_{2}\neq 0$,
then
\begin{eqnarray*}
&\ &\textrm{ord}_{p}\omega_{1}=\frac{p+7}{6(p-1)},
\textrm{ord}_{p}\omega_{2}=\frac{p-2}{3(p-1)},
\textrm{ord}_{p}\omega_{3}=\frac{1}{2},
\textrm{ord}_{p}\omega_{4}=1-\frac{p-2}{3(p-1)},\\
&\ &\textrm{ord}_{p}\omega_{5}=1-\frac{p+7}{6(p-1)}.
\end{eqnarray*}

\textrm{(vii)} If $a_{4}=a_{3}= 0$, $a_{2}\neq 0$, then
\begin{eqnarray*}
&\ &\textrm{ord}_{p}\omega_{1}=\frac{p+7}{6(p-1)},
\textrm{ord}_{p}\omega_{2}=\frac{p+1}{3(p-1)},
\textrm{ord}_{p}\omega_{3}=\frac{1}{2},
\textrm{ord}_{p}\omega_{4}=1-\frac{p+1}{3(p-1)},\\
&\ &\textrm{ord}_{p}\omega_{5}=1-\frac{p+7}{6(p-1)}.
\end{eqnarray*}

\textrm{(viii)} If $a_{4}=0$, $a_{3}\neq 0$,
$a_{2}^{2}+2a_{1}a_{3}=0$, $9a_{2}^{3}-10a_{3}^{4}\neq 0$, and
$a_{1}\neq 0, a_{2}\neq 0$, then
\begin{eqnarray*}
&\ &\textrm{ord}_{p}\omega_{1}=\frac{p+13}{6(p-1)},
\textrm{ord}_{p}\omega_{2}=\frac{p-5}{3(p-1)},
\textrm{ord}_{p}\omega_{3}=\frac{1}{2},
\textrm{ord}_{p}\omega_{4}=1-\frac{p-5}{3(p-1)},\\
&\ &\textrm{ord}_{p}\omega_{5}=1-\frac{p+13}{6(p-1)}.
\end{eqnarray*}

\textrm{(ix)} If $a_{4}=0$, $a_{3}\neq 0$,
$a_{2}^{2}+2a_{1}a_{3}=0$, $9a_{2}^{3}-10a_{3}^{4}=0$, then
\begin{eqnarray*}
&\ &\textrm{ord}_{p}\omega_{1}=\frac{p+19}{6(p-1)},
\textrm{ord}_{p}\omega_{2}=\frac{p-8}{3(p-1)},
\textrm{ord}_{p}\omega_{3}=\frac{1}{2},
\textrm{ord}_{p}\omega_{4}=1-\frac{p-8}{3(p-1)},\\
&\ &\textrm{ord}_{p}\omega_{5}=1-\frac{p+19}{6(p-1)}.
\end{eqnarray*}

\textrm{(x)} If $a_{4}=a_{3}=a_{2}=0$, $a_{1}\neq 0$, then
\begin{eqnarray*}
&\ &\textrm{ord}_{p}\omega_{1}=\frac{p+19}{6(p-1)},
\textrm{ord}_{p}\omega_{2}=\frac{p+4}{3(p-1)},
\textrm{ord}_{p}\omega_{3}=\frac{1}{2},
\textrm{ord}_{p}\omega_{4}=1-\frac{p+4}{3(p-1)},\\
&\ &\textrm{ord}_{p}\omega_{5}=1-\frac{p+19}{6(p-1)}.
\end{eqnarray*}

\textrm{(xi)} If $a_{4}=a_{2}=a_{1}=0$, $a_{3}\neq 0$, then
$$
\textrm{ord}_{p}\omega_{1}=\textrm{ord}_{p}\omega_{2}=\frac{p+1}{4(p-1)},
\textrm{ord}_{p}\omega_{3}=\frac{1}{2},
\textrm{ord}_{p}\omega_{4}=\textrm{ord}_{p}\omega_{5}=1-\frac{p+1}{4(p-1)}.
$$
\end{Theorem}

\begin{pf}
Similar as above, we can generally obtain that
\begin{eqnarray*}
C(\frac{p+1}{6}; 1, 1)=\{[0, 0, 0, 1]\},
\end{eqnarray*}

\begin{eqnarray*}
C(1+\frac{p+1}{6}; 1, 1)=\{[0, 0, 0, 4], [0, 1, 0, 2], [0, 2, 0, 0],
[1, 0, 1, 0], [0, 0, 2, 1]\},
\end{eqnarray*}

\begin{eqnarray*}
C(2+\frac{p+1}{6}; 1, 1)&=&\{[0, 0, 0, 7], [2, 0, 0, 2], [0, 1, 0,
5], [2, 1, 0, 0], [0, 2, 0, 3], [0, 3, 0, 1], \\
&\ &[1, 0, 1, 3], [1, 1, 1, 1], [0, 0, 2, 4], [0, 1, 2, 2], [0, 2,
2, 0], [1, 0, 3, 0], \\
&\ &[0, 0, 4, 1]\},
\end{eqnarray*}

\begin{eqnarray*}
C(3+\frac{p+1}{6}; 1, 1)&=&\{[0, 0, 0, 10], [2, 0, 0, 5], [4, 0, 0,
0], [0, 1, 0, 8], [2, 1, 0, 3], [0, 2, 0, 6], \\
&\ &[2, 2, 0, 1], [0, 3, 0, 4], [0, 4, 0, 2], [0, 5, 0, 0], [1, 0,
1, 6], [3, 0, 1, 1], \\
&\ &[1, 1, 1, 4], [1, 2, 1, 2], [1, 3, 1, 0], [0, 0, 2, 7], [2, 0,
2, 2], [0, 1, 2, 5], \\
&\ &[2, 1, 2, 0], [0, 2, 2, 3], [0, 3, 2, 1], [1, 0, 3, 3], [1, 1,
3, 1], [0, 0, 4, 4], \\
&\ &[0, 1, 4, 2], [0, 2, 4, 0], [1, 0, 5, 0], [0, 0, 6, 1]\},
\end{eqnarray*}
and
\begin{eqnarray*}
C(\frac{p+1}{3}; 2, 2)=\{[0, 0, 0, 2], [0, 1, 0, 0]\},
\end{eqnarray*}

\begin{eqnarray*}
C(1+\frac{p+1}{3}; 2, 2)&=&\{[0, 0, 0, 5], [2, 0, 0, 0], [0, 1, 0,
3], [0, 2, 0, 1], [1, 0, 1, 1], [0, 0, 2, 2], \\
&\ &[0, 1, 2, 0]\},
\end{eqnarray*}

\begin{eqnarray*}
C(2+\frac{p+1}{3}; 2, 2)&=&\{[0, 0, 0, 8], [2, 0, 0, 3], [0, 1, 0,
6], [2, 1, 0, 1], [0, 2, 0, 4], [0, 3, 0, 2], \\
&\ &[0, 4, 0, 0], [1, 0, 1, 4], [1, 1, 1, 2], [1, 2, 1, 0], [0, 0,
2, 5], [2, 0, 2, 0], \\
&\ &[0, 1, 2, 3], [0, 2, 2, 1], [1, 0, 3, 1], [0, 0, 4, 2], [0, 1,
4, 0]\},
\end{eqnarray*}

\begin{eqnarray*}
C(3+\frac{p+1}{3}; 2, 2)&=&\{[0, 0, 0, 11], [2, 0, 0, 6], [4, 0, 0,
1], [0, 1, 0, 9], [2, 1, 0, 4], [0, 2, 0, 7], \\
&\ &[2, 2, 0, 2], [0, 3, 0, 5], [2, 3, 0, 0], [0, 4, 0, 3], [0, 5,
0, 1], [1, 0, 1, 7], \\
&\ &[3, 0, 1, 2], [1, 1, 1, 5], [3, 1, 1, 0], [1, 2, 1, 3], [1, 3,
1, 1], [0, 0, 2, 8], \\
&\ &[2, 0, 2, 3], [0, 1, 2, 6], [2, 1, 2, 1], [0, 2, 2, 4], [0, 3,
2, 2], [0, 4, 2, 0], \\
&\ &[1, 0, 3, 4], [1, 1, 3, 2], [1, 2, 3, 0], [0, 0, 4, 5], [2, 0,
4, 0], [0, 1, 4, 3], \\
&\ &[0, 2, 4, 1], [1, 0, 5, 1], [0, 0, 6, 2], [0, 1, 6, 0]\},
\end{eqnarray*}
and
\begin{eqnarray*}
C(\frac{p+1}{6}; 1, 2)=\{[0, 0, 1, 0]\},
\end{eqnarray*}

\begin{eqnarray*}
C(1+\frac{p+1}{6}; 1, 2)=\{[1, 0, 0, 2], [1, 1, 0, 0], [0, 0, 1, 3],
[0, 1, 1, 1], [0, 0, 3, 0]\},
\end{eqnarray*}

\begin{eqnarray*}
C(2+\frac{p+1}{6}; 1, 2)&=&\{[1, 0, 0, 5], [3, 0, 0, 0], [1, 1, 0,
3], [1, 2, 0, 1], [0, 0, 1, 6], [2, 0, 1, 1], \\
&\ &[0, 1, 1, 4], [0, 2, 1, 2], [0, 3, 1, 0], [1, 0, 2, 2], [1, 1,
2, 0], [0, 0, 3, 3], \\
&\ &[0, 1, 3, 1], [0, 0, 5, 0]\},
\end{eqnarray*}

\begin{eqnarray*}
C(3+\frac{p+1}{6}; 1, 2)&=&\{[1, 0, 0, 8], [3, 0, 0, 3], [1, 1, 0,
6], [3, 1, 0, 1], [1, 2, 0, 4], [1, 3, 0, 2], \\
&\ &[1, 4, 0, 0], [0, 0, 1, 9], [2, 0, 1, 4], [0, 1, 1, 7], [2, 1,
1, 2], [0, 2, 1, 5], \\
&\ &[2, 2, 1, 0], [0, 3, 1, 3], [0, 4, 1, 1], [1, 0, 2, 5], [3, 0,
2, 0], [1, 1, 2, 3], \\
&\ &[1, 2, 2, 1], [0, 0, 3, 6], [2, 0, 3, 1], [0, 1, 3, 4], [0, 2,
3, 2], [0, 3, 3, 0], \\
&\ &[1, 0, 4, 2], [1, 1, 4, 0], [0, 0, 5, 3], [0, 1, 5, 1], [0, 0,
7, 0]\},
\end{eqnarray*}
and
\begin{eqnarray*}
C(\frac{p+1}{3}; 2, 1)=\{[0, 0, 1, 0]\},
\end{eqnarray*}

\begin{eqnarray*}
C(1+\frac{p+1}{3}; 2, 1)=\{[1, 0, 0, 2], [1, 1, 0, 0], [0, 0, 1, 3],
[0, 1, 1, 1], [0, 0, 3, 0]\},
\end{eqnarray*}

\begin{eqnarray*}
C(2+\frac{p+1}{3}; 2, 1)&=&\{[1, 0, 0, 5], [3, 0, 0, 0], [1, 1, 0,
3], [1, 2, 0, 1], [0, 0, 1, 6], [2, 0, 1, 1], \\
&\ &[0, 1, 1, 4], [0, 2, 1, 2], [0, 3, 1, 0], [1, 0, 2, 2], [1, 1,
2, 0], [0, 0, 3, 3], \\
&\ &[0, 1, 3, 1], [0, 0, 5, 0]\},
\end{eqnarray*}

\begin{eqnarray*}
C(3+\frac{p+1}{3}; 2, 1)&=&\{[1, 0, 0, 8], [3, 0, 0, 3], [1, 1, 0,
6], [3, 1, 0, 1], [1, 2, 0, 4], [1, 3, 0, 2], \\
&\ &[1, 4, 0, 0], [0, 0, 1, 9], [2, 0, 1, 4], [0, 1, 1, 7], [2, 1,
1, 2], [0, 2, 1, 5], \\
&\ &[2, 2, 1, 0], [0, 3, 1, 3], [0, 4, 1, 1], [1, 0, 2, 5], [3, 0,
2, 0], [1, 1, 2, 3], \\
&\ &[1, 2, 2, 1], [0, 0, 3, 6], [2, 0, 3, 1], [0, 1, 3, 4], [0, 2,
3, 2], [0, 3, 3, 0], \\
&\ &[1, 0, 4, 2], [1, 1, 4, 0], [0, 0, 5, 3], [0, 1, 5, 1], [0, 0,
7, 0]\}.
\end{eqnarray*}

Thus we can calculate every $F_{r}^{s}$ ($s=1, 2$) and get
conclusions by discussion on whether $F_{r}^{s}=0$ for some $r$ and
$s$.
\end{pf}

\begin{Remark}
When $p\leq 29$, some $C(r; u, v)$ in the proof will lost some
elements. For example, when $p=23$, $C(2+\frac{p+1}{6}; 1, 1)$ will
lost $[0, 0, 0, 7]$ and $C(3+\frac{p+1}{6}; 1, 1)$ will lost $[0, 0,
0, 10], [0, 1, 0, 8], [0, 2, 0, 6], [1, 0, 1, 6], [0, 0, 2, 7], [0,
0, 4, 4]$. The reason is, $r-\sum_{i=1}^{m-1}k_{i}$ must be a
non-negative integer. So these cases when $p$ are small should be
considered specially, see case $(ii)$ and $(iii)$.

Most cases in \textbf{theorem 7.5.} are consistent with S.Hong's
result (see [10]), but some cases did not discuss by S.Hong. For
example, S.Hong lost the case $(ii)$, it is possible since $a_{1}=2,
a_{2}=6, a_{3}=a_{4}=3$ is such an example.

Furthermore the first case in S.Hong's Theory
($a_{1}a_{2}a_{3}a_{4}\neq 0$) had some mistake for the same reason.
\end{Remark}

\section{\bf Some other examples}

To illustrate our method, consider the case where $f=x^{7}+ax^{4}$,
$a\in \textbf{F}_{p}^{*}$ and $p=5$ a small prime. Similar to the
method in \textbf{section 3}, we also have a table to indicate the
factor in (7), the only different is, it do not need satisfy the
relation (9).

It is clear that $L^{*}(f/\textbf{F}_{p},T)$ is a polynomial of
degree 7. We need to consider the solutions $r=(r_{1},r_{2})\in
S_{p}(f)$ of equation
\begin{equation}
\left(
\begin{array}{cccc}
4, 7\\
\end{array}
\right)\left(
\begin{array}{cccc}
r_{1}\\
r_{2}\\
\end{array}
\right)\equiv 0 \ \textrm{mod}\  1.
\end{equation}

Consider the coefficient $c_{2}$, $\left(
\begin{array}{cccc}
0\\
0\\
\end{array}
\right)\left(
\begin{array}{cccc}
1\\
0\\
\end{array}
\right)$ is the only table which satisfies (17) and makes the sum of
all elements in the table smallest. So we have
$$
\textrm{ord}_{p}c_{2}=\frac{1}{4}.
$$

Consider the coefficient $c_{3}$, $\left(
\begin{array}{cccc}
0\\
0\\
\end{array}
\right)\left(
\begin{array}{cccc}
1\\
0\\
\end{array}
\right)
\left(
\begin{array}{cccc}
2\\
0\\
\end{array}
\right)$ is the only table which satisfies (17) and makes the sum of
all elements in the table smallest. So we have
$$
\textrm{ord}_{p}c_{3}=\frac{3}{4}.
$$

Consider the coefficient $c_{4}$, $\left(
\begin{array}{cccc}
0\\
0\\
\end{array}
\right)\left(
\begin{array}{cccc}
1\\
0\\
\end{array}
\right)\left(
\begin{array}{cccc}
2\\
0\\
\end{array}
\right)\left(
\begin{array}{cccc}
3\\
0\\
\end{array}
\right)$ and $\left(
\begin{array}{cccc}
0\\
0\\
\end{array}
\right)\left(
\begin{array}{cccc}
1\\
0\\
\end{array}
\right)\left(
\begin{array}{cccc}
1,0\\
0,4\\
\end{array}
\right)$ and $\left(
\begin{array}{cccc}
0\\
0\\
\end{array}
\right)\left(
\begin{array}{cccc}
0,1,1\\
1,1,2\\
\end{array}
\right)$ and $\left(
\begin{array}{cccc}
0\\
0\\
\end{array}
\right)\left(
\begin{array}{cccc}
0,2,0\\
0,2,2\\
\end{array}
\right)$ are the only 4 tables which satisfies (17) and makes the
sum of all elements in the table smallest respectively. Then the sum
of their $\Gamma_{p}$-parts $\ \textrm{mod}\  p$ is equal to
$$
\frac{1}{2!}\cdot\frac{1}{3!}\cdot\chi(a)^{1+2+3}-\frac{1}{4!}\cdot\chi(a)^{1+1}+\frac{1}{2!}\cdot\chi(a)^{1+1}+\frac{1}{2!}\cdot\frac{1}{2!}\cdot\frac{1}{2!}\cdot\chi(a)^{2},
$$
that is
$$
\frac{1}{12}\chi(a)^{2}(\chi(a)^{4}+7).
$$

Since $a\in \textbf{F}_{5}^{*}$, we have $\chi(a)^{4}=1$ and
therefore $\frac{1}{12}\chi(a)^{2}(\chi(a)^{4}+7)\neq 0 \
\textrm{mod}\ 5$. That means
$$
\textrm{ord}_{p}c_{4}=\frac{3}{2}.
$$

For the reason that the Newton polygon of $L(f,T)$ is symmetric in
the sense that for every slope segment $\alpha$ there is a slope
segment $1-\alpha$ of the same horizontal length, and
$\textrm{ord}_{p}c_{3}-\textrm{ord}_{p}c_{2}=\frac{1}{2}$,
$\textrm{ord}_{p}c_{4}-\textrm{ord}_{p}c_{3}=\frac{3}{4}>\frac{1}{2}$,
we have
$$
\textrm{ord}_{p}\omega_{1}=\frac{1}{4},
$$
$$
\textrm{ord}_{p}\omega_{2}=\textrm{ord}_{p}\omega_{3}=\textrm{ord}_{p}\omega_{4}=\textrm{ord}_{p}\omega_{5}=\frac{1}{2},
$$
$$
\textrm{ord}_{p}\omega_{6}=\frac{3}{4}.
$$
\begin{center}
\end{center}

\begin{center}
\end{center}

An other example is for $f=x^{3}+axy+by^{2}$, where $q=p^{k}$ and
$a, b\in \textbf{F}_{q}^{*}$ and $p>6$ is prime satisfying $p\equiv
-1 \ \textrm{mod}\  3$.

It is clear that $L^{*}(f/\textbf{F}_{q},T)^{(-1)}$ is a polynomial
of degree 6. We need to consider the solutions
$r=(r_{1},r_{2},r_{3})\in S_{p}(f)$ of equation
\begin{equation}
\left(
\begin{array}{cccc}
3, 1, 0\\
0, 1, 2\\
\end{array}
\right)\left(
\begin{array}{cccc}
r_{1}\\
r_{2}\\
r_{3}\\
\end{array}
\right)\equiv 0 \ \textrm{mod}\  1.
\end{equation}

Consider the coefficient $c_{2}$, the solutions $ \left(
\begin{array}{cccc}
0\\
0\\
0\\
\end{array}
\right)$ and $\left(
\begin{array}{cccc}
0\\
0\\
\frac{1}{2}\\
\end{array}
\right)$ of (18) are made up of a table below:

$$
\left(
\begin{array}{cccc}
0, \cdots, 0\\
0, \cdots, 0\\
0, \cdots, 0\\
\end{array}
\right)
\left(
\begin{array}{cccc}
0, \cdots, 0\\
0, \cdots, 0\\
\frac{p-1}{2}, \cdots, \frac{p-1}{2}\\
\end{array}
\right)
$$
Each block in this table has $k$ columns. It is the only table which
satisfies (18) and makes the sum of all elements in it smallest. So
we have
$$
\textrm{ord}_{q}c_{2}=\frac{1}{2}.
$$

Consider the coefficient $c_{3}$, the solutions $ \left(
\begin{array}{cccc}
0\\
0\\
0\\
\end{array}
\right)$ and $\left(
\begin{array}{cccc}
\frac{1}{3}\\
0\\
0\\
\end{array}
\right)$ of (18) are made up of a table below:

(i) $k\equiv 1 \ \textrm{mod}\  2$, then
$$
\left(
\begin{array}{cccc}
0, \cdots, 0\\
0, \cdots, 0\\
0, \cdots, 0\\
\end{array}
\right)\left(
\begin{array}{cccc}
\frac{2p-1}{3}, \frac{p-2}{3}, \cdots, \frac{2p-1}{3}, \frac{p-2}{3}\\
0, \ \ \ 0, \ \ \  \cdots, \ \ \ 0, \ \ \ 0\\
0, \ \ \ 0, \ \ \  \cdots, \ \ \ 0, \ \ \ 0\\
\end{array}
\right)
$$

The first block in this table has $k$ columns and the second block
has $2k$ columns.

(ii) $k\equiv 0 \ \textrm{mod}\  2$, then
$$
\left(
\begin{array}{cccc}
0, \cdots, 0\\
0, \cdots, 0\\
0, \cdots, 0\\
\end{array}
\right)\left(
\begin{array}{cccc}
\frac{p-2}{3}, \frac{2p-1}{3}, \cdots, \frac{p-2}{3}, \frac{2p-1}{3}\\
0, \ \ \ 0, \ \ \  \cdots, \ \ \ 0, \ \ \ 0\\
0, \ \ \ 0, \ \ \  \cdots, \ \ \ 0, \ \ \ 0\\
\end{array}
\right)\left(
\begin{array}{cccc}
\frac{2p-1}{3}, \frac{p-2}{3}, \cdots, \frac{2p-1}{3}, \frac{p-2}{3}\\
0, \ \ \ 0, \ \ \  \cdots, \ \ \ 0, \ \ \ 0\\
0, \ \ \ 0, \ \ \  \cdots, \ \ \ 0, \ \ \ 0\\
\end{array}
\right)
$$
Each block in this table has $k$ columns.

It is the only table which satisfies (18) and makes the sum of all
elements in it smallest. So we have
$$
\textrm{ord}_{q}c_{3}=1.
$$

Consider the coefficient $c_{4}$, we can easily prove that add
$\left(
\begin{array}{cccc}
0, \cdots, 0\\
0, \cdots, 0\\
\frac{p-1}{2}, \cdots, \frac{p-1}{2}\\
\end{array}
\right) $ to the table (i) or (ii) is the only table which satisfies
(18) and makes the sum of all elements in it smallest. So we have
$$
\textrm{ord}_{q}c_{4}=\frac{3}{2}.
$$

Consider the coefficient $c_{5}$ and note the solution $ \left(
\begin{array}{cccc}
\frac{p-5}{3(p-1)}\\
\frac{4}{p-1}\\
\frac{p-5}{2(p-1)}\\
\end{array}
\right) $of (18). We can easily prove that add $\left(
\begin{array}{cccc}
\frac{p-5}{3}, \cdots, \frac{p-5}{3}\\
4, \cdots, 4\\
\frac{p-5}{2}, \cdots, \frac{p-5}{2}\\
\end{array}
\right) $ to the table of case $c_{4}$ is the only table which
satisfies (18) and makes the sum of all elements in it smallest. So
we have
$$
\textrm{ord}_{q}c_{5}=\frac{7}{3}+\frac{2}{3(p-1)}.
$$

We can not calculate $\textrm{ord}_{q}c_{6}$ by our method yet but
$\textrm{ord}_{p}c_{6}$ since

$$\left(
\begin{array}{cccc}
0\\
0\\
0\\
\end{array}
\right)\left(
\begin{array}{cccc}
0\\
0\\
\frac{p-1}{2}\\
\end{array}
\right)\left(
\begin{array}{cccc}
\frac{2p-1}{3}, \frac{p-2}{3}\\
0, 0\\
0, 0\\
\end{array}
\right)\left(
\begin{array}{cccc}
\frac{p-5}{3}\\
4\\
\frac{p-5}{2}\\
\end{array}
\right)\left(
\begin{array}{cccc}
\frac{2p-4}{3}\\
2\\
\frac{p-3}{2}\\
\end{array}
\right)
$$

is the only table which satisfies (18) and makes the sum of all
elements in it smallest. So we have
$$
\textrm{ord}_{p}c_{6}=\frac{7}{2}+\frac{1}{p-1}.
$$

Then for $\frac{1}{1-T}\cdot L^{*}(f/\textbf{F}_{q},T)^{(-1)}$ we
have
$$
\textrm{ord}_{q}\omega_{1}=\textrm{ord}_{q}\omega_{2}=\textrm{ord}_{p}\omega_{3}=\frac{1}{2},
$$
$$
\textrm{ord}_{q}\omega_{4}=\frac{5}{6}+\frac{2}{3(p-1)}.
$$
and for $q=p$ we have
$$
\textrm{ord}_{p}\omega_{5}=\frac{7}{6}+\frac{1}{3(p-1)}.
$$

\begin{Remark}
our calculations here are based on  \textbf{theorem 2.2.} and
\textbf{section 3}. We can also see that the table-representation
and \textbf{proposition 3.4.} are applicable not only for $q=p$ or
$n=1$.
\end{Remark}
\begin{center}
\end{center}

\textbf{Acknowledgements} The authors thank professor Daqing Wan
deeply for introducing them into this fascinating field by offering
excellent lecture notes of his, answering many questions in their
seminar or by emails for many times.

\medskip

\end{document}